\documentclass[journal]{IEEEtran}
\usepackage{amssymb,times}
\usepackage{lipsum}%
\usepackage{cite, epsfig}
\usepackage{verbatim, amsfonts,enumerate}     
\usepackage{bm,array}
\usepackage{graphicx}
\usepackage[cmex10]{amsmath}
\usepackage{epstopdf}
\usepackage{subcaption}
\usepackage{amsthm}
\usepackage{algorithm,algorithmic}
\usepackage{framed}
\usepackage{color,url}
\usepackage[normalem]{ulem}
\usepackage{graphicx}

\theoremstyle{plain}
\newtheorem{thm}{Theorem}
\newtheorem{lem}{Lemma}
\DeclareMathOperator*{\argmin}{arg\,min}
\def \S {\mathcal{S}}

\def\v{\mathbf{v}}
\def\u{\mathbf{u}}

\def\C{\mathbf{C}}
\def\X{\mathcal{X}}
\def\x{\mathbf{x}}

\def\Q{\mathbf{Q}}

\def\A{\mathbf{A}}
\def\fh{\hat{f}}
\def\gh{\hat{g}}
\def\y{\mathbf{y}}
\def\e{\mathbf{e}}
\def\yh{\hat{\mathbf{y}}}

\def\R{\mathbb{R}}
\def\z{\mathbf{z}}

\def\xh{\hat\x}
\def\zh{\hat{\z}}

\def\prox{\text{prox}}
\def\G{\mathcal{G}}
\def\T{\mathcal{T}}
\def\N{\mathcal{N}}
\def\E{\mathcal{E}}
\def\xb{\bar\xh}
\def\zb{\bar{\zh}}
\def\ve{\boldsymbol{\varepsilon}}
\def\Reg{\mathbf{Reg}_T^D}
\def\t{\lfloor\! t\!\rfloor}

\providecommand{\abs}[1]{\left|#1\right|}
\providecommand{\norm}[1]{\left\|#1\right\|}
\providecommand{\ip}[1]{\langle#1\rangle}

\newcommand{\colb}[1]{\textcolor{black}{#1}}

\theoremstyle{plain}
\theoremstyle{remark}

\newtheorem{assum}{\bf Assumption\!}  

\newtheorem{cor}{\bf Corollary}

\graphicspath{{img/}}
\begin{document}

\title{Online Learning over Dynamic Graphs via Distributed Proximal Gradient Algorithm}
\author{Rishabh Dixit,	Amrit Singh Bedi,~\IEEEmembership{Student Member,~IEEE},  
	and~Ketan~Rajawat,~\IEEEmembership{Member,~IEEE}
	\thanks{
		The authors are with the Department of Electrical Engineering,
		Indian Institute of Technology Kanpur, Kanpur 208016, India (e-mail:
		\texttt{rishabd93@gmail.com}; \texttt{amritbd@iitk.ac.in}; \texttt{ketan@iitk.ac.in}). This work was supported by the DST-SERB Grant EMR/2016/005959.}\vspace{-0mm}}
	\maketitle
\begin{abstract}
We consider the problem of tracking the minimum of a time-varying convex optimization problem over a dynamic graph. Motivated by target tracking and parameter estimation problems in intermittently connected robotic and sensor networks, the goal is to design a distributed algorithm capable of handling non-differentiable regularization penalties. The proposed proximal online gradient descent algorithm is built to run in a fully decentralized manner and utilizes consensus updates over possibly disconnected graphs. The performance of the proposed algorithm is analyzed by developing bounds on its dynamic regret in terms of the cumulative path length of the time-varying optimum. It is shown that as compared to the centralized case, the dynamic regret incurred by the proposed algorithm over $T$ time slots is worse by a factor of $\log(T)$ only, despite the disconnected and time-varying network topology. The empirical performance of the proposed algorithm is tested on the distributed dynamic sparse recovery problem, where it is shown to incur a dynamic regret that is close to that of the centralized algorithm. 
\end{abstract}

\begin{IEEEkeywords} Dynamic regret, distributed optimization, online convex optimization, sparse signal recovery. 
\end{IEEEkeywords}

\section{Introduction}
Recent advances in sensing, communication, and computation technologies have ushered an era of large-scale networked systems, now ubiquitous in smart grids, robotics, defense installations, industrial automation, and cyber-physical systems \cite{chen2010distributed, dall2013distributed, chen2017online}. In order to accomplish sophisticated tasks, these multi-agent systems  are expected to sense, learn, and adapt to sequentially arriving measurements while facing  time-varying and non-stationary environments. Such high levels of network agility requires close coordination and cooperation among the agents, that are otherwise intermittently connected and resource-constrained. Towards performing control and resource allocation over such dynamic graphs, it is necessary to develop in-network optimization algorithms that can be implemented in a distributed fashion \cite{nedic2018distributed}. 

We consider the problem of tracking the minimum of a time-varying cost function that can be written as a sum of several node-specific smooth cost functions and a non-smooth regularizer. The framework subsumes a number of relevant problems such as model-predictive control \cite{jerez2014embedded, gutjahr2017lateral}, tracking time-varying parameters \cite{jakubiec2013d, vaswani2016recursive, simonetto2015non}, path planning \cite{ardeshiri2011convex, verscheure2009time}, real-time magnetic resonance imaging \cite{uecker2012real}, adaptive matrix completion \cite{tripathi2017adaptive}, demand scheduling in smart grids \cite{zhao2014design}, and so on. Of particular interest is the distributed setting, where measurements are made locally at each node and the interaction between nodes may only occur over a time-varying directed communication graph. Due to the intermittent connectivity, it is not possible for the nodes to exchange updates at every time slot, and consequently, simple variants of centralized algorithms cannot be directly utilized.

Within the context of static optimization problems, distributed algorithms running over time-varying graphs have been widely studied, with most algorithms relying on the consensus approach \cite{nedich2015convergence}. Consensus-based algorithms can be further categorized into those utilizing weighted averaging in the primal domain \cite{nedic2009distributed,nedic2015distributed,ram2009incremental}, push-sum-based approaches \cite{benezit2010weighted}, and the alternating directions method of multipliers (ADMM) method \cite{boyd2011distributed}. The weighted-averaging-based approaches utilize a doubly stochastic matrix for averaging and are amenable to gradient-tracking variants that converge at geometric rates; see \cite{nedic2018distributed} and references therein. 

Fewer distributed algorithms exist for the more challenging time-varying setting where the minimizer drifts over time. Algorithms for tracking time-varying parameters have their roots in adaptive signal processing and control theory, where the steady-state tracking error is of interest \cite{widrow1976stationary,7902101}. Recently however, online convex optimization has emerged as a useful framework for the analysis of tracking algorithms, especially in adversarial settings \cite{zinkevich2003online, hall2015online, besbes2015non, jadbabaie2015online, mokhtari2016online, yang2016tracking, zhang2016improved, bedi2018tracking}. The idea here is to characterize the dynamic regret of an algorithm in terms of the cumulative variations in the problem parameters such as the path length. However, existing online algorithms for tracking time-varying objectives cannot handle time-varying graph topologies where the nodes may get disconnected. A distributed dynamic ADMM algorithm for differentiable cost functions was first proposed in \cite{ling2014decentralized} and the steady-state optimality gap was derived for static graphs. On the other hand, an online ADMM algorithm for  non-differentiable problems was proposed in \cite{simonetto2015non} but required centralized implementation. Closer to the current work, the dynamic regret performance of a distributed mirror descent algorithm was obtained in \cite{shahrampour2018distributed} for possibly non-differentiable functions. However, the underlying graph was required to be static and connected. In contrast, the focus here is on distributed algorithms over dynamic and sporadically connected graphs. 

This works puts forth a distributed proximal online gradient descent (OGD) algorithm designed to run over a time-varying network topology and capable of tracking the minimum of a time-varying composite loss function. As the objective may contain non-differentiable regularizers, we build upon the recently developed machinery of proximal OGD \cite{dixit2018online}. Towards realizing a distributed implementation, we utilize the idea of weighted averaging using a doubly stochastic weight matrix, as is customary in distributed algorithms for static problems \cite{nedic2018distributed}. However, since the graph topology is time-varying, it is not possible for the nodes to perform a consensus update or even exchange information at each time slot. Instead a multi-step consensus algorithm inspired from \cite{chen2012fast} is developed, where the objective function information is used only intermittently and the remaining time slots are utilized for consensus. It is remarked that the time-varying problem is significantly more challenging than the static version considered in \cite{chen2012fast}, since the objective function continues to evolve even when the consensus steps are being carried out, resulting in the accumulation of additional regret. 

The proposed distributed proximal OGD algorithm is analyzed as an inexact variant of the proximal OGD algorithm \cite{rates11schmidt}. Subsequently, the gradient and proximal errors arising due to the distributed operation are characterized, allowing us to obtain the required dynamic regret bounds. Interestingly, the distributed operation only results in an additional factor of $\log(T)$ as compared to the dynamic regret lower bound for the centralized case. The performance of the proposed algorithm is tested for the dynamic sparse recovery problem and compared with that of the centralized proximal OGD and ADMM algorithms. 

\colb{A large number of LMS- and RLS-based algorithm have been proposed to solve the dynamic sparse recovery problem, starting from \cite{angelosante2010online,YilunChen:2009:SLS:1582709.1583473}. An online ADMM approach for time-varying LASSO problem was developed in  \cite{wang2012online} while a online proximal ADMM algorithm for solving the group LASSO problem was proposed in \cite{suzuki2013dual}. While analytical results have been developed, they are largely limited to the case of static parameters. In particular, these works do not generally characterize the steady state tracking performance, though it may be possible to borrow similar bounds from the matrix completion literature; see  \cite{tripathi2017adaptive}. A comprehensive survey of RLS-based dynamic sparse recovery methods is provided in \cite{vaswani2016recursive}. }


\colb{When the sparse parameters follow a linear state-space model, they can be tracked within the sparse Bayesian learning (SBL) formalism \cite{6638927,6288632}. More recently, consensus has been employed to track the sparse parameters in a distributed fashion \cite{khanna2017decentralized}. Likewise, a robust framework for dynamic sparse recovery has been developed in \cite{radhika2018robust}. Different from the system model considered here, SBL-based methods cannot handle adversarial targets and are generally not amenable to regret analysis. }


\colb{Adversarial parameters are naturally handled within the dynamic or online convex optimization framework. For instance, the dynamic sparse recovery problem was formulated and solved via a 'running' ADMM approach in \cite{simonetto2015non}. Likewise, the mirror descent algorithm for solving a similar problem was proposed in \cite{shahrampour2018distributed}. Generalizing these frameworks, we consider the dynamic sparse recovery problem in a  distributed setting where the number of observations per sensor are not sufficient to estimate the full parameter. Further, the underlying graph topology is time-varying and collecting the observations at a central location is not straightforward. }

In summary, our key contributions include: (a) a novel multi-step consensus-based proximal OGD algorithm for time-varying optimization problems over time-varying network topologies (b) performance analysis of the general inexact proximal OGD algorithm and its application to the analysis of the proposed distributed algorithm. The rest of the paper is organized as follows. We start with the problem formulation in Sec. \ref{Prob_for}. The proposed distributed proximal OGD algorithm is presented in Sec. \ref{seciii}. A bound on the dynamic regret incurred by the proposed algorithm is obtained in Sec. \ref{perfa}. The empirical performance of the proposed algorithm is tested on the dynamic sparse recovery problem in Sec. \ref{numerical} and the conclusions are presented in Sec.\ref{conclusion}. 

{Notation:} All the scalars are denoted by regular font lower case letters, vectors by boldface lower case letters, and matrices by upper case boldface letters. Constants such as the total number of time slots and gradient bounds are denoted by regular font capital letters (such as $T$ and $M$) while problem parameters such as step-size are denoted by Greek letters (such as $\alpha$). Finally, sets are denoted by calligraphic font capital letters.  The  $(i,j)$-th entry of the matrix $\A \in \R^{m \times n}$ is denoted either by $[\A]_{ij}$  or $A^{ij}$. The Euclidean norm of a vector $\x \in \R^n$ is denoted by $\norm{\cdot}$. By default, time and iteration indices (e.g. $t$ or $k$) will be in the subscript, while node or element indices (e.g. $i$ or $j$) will be in the superscript.

\begin{table*}
			\centering
			\captionof{table}{Dynamic regret rates for online learning with non-differentiable functions (cf. Sec.\ref{seciii})}	\label{table} \hspace*{-0.3cm}
			\setlength\tabcolsep{1.5pt}
			\begin{tabular}{cccccc}
				\hline
				References &  Convex & Distributed & Dynamic graph topology & Dynamic regret bound \\
				\hline 
				\cite{hall2015online} 		& Convex				& No & No &	$\mathcal{O}(\sqrt{T}(1+C_T)$  \\ 
				\cite{shahrampour2018distributed} 					& Convex & Yes	& No		 	& 	 ${\mathcal{O}{(\sqrt{TC_T})}}$ \\
				This work 					& Strongly convex  			 	& Yes & Yes   & $\mathcal{O}(\log{(T)}(1+C_T))$  \\
	\hline
	\end{tabular}
\end{table*}

%

	\section{Problem Formulation}\label{Prob_for}
Consider a multi-agent network with the set of agents or nodes $\N:=\{1, \ldots, N\}$ interacting with each other over intermittent communication links. The network topology at time $t$ is represented by a directed graph $\G_t :=(\N, \E_t)$ where $\E_t \subseteq \N \times \N$ denotes the set of directed edges or links present at time $t$. Specifically, agents $i$ can transmit to agent $j$ only if $(i,j) \in \E_t$ at time instant $t$. The dynamic graph $\G_t$ is arbitrary and may not necessarily be connected for all $t\geq 1$. However, $\{\G_t\}$ is still required to be uniformly strongly connected or $B$-connected \cite{nedic2015distributed}; see Sec. \ref{seciii}. 
	
Restricted to exchanging information only over the edges in $\E_t$ at time $t$, the agents seek to cooperatively track the optimum of the following dynamic optimization problem
\begin{align}\label{p}
\x_t^\star:=\arg\min_{\x\in\X_t} \ell_t(\x):=\frac{1}{N}\sum_{i=1}^N f_t^{i}(\x) + h_t(\x)\tag{$\mathcal{P}$}
\end{align}
where $f_t^i:\R^n\rightarrow \R$ is a smooth strongly convex function,  $h_t: \R^n \rightarrow \R$ is a convex but possibly non-smooth regularizer, and $\X_t \subset \R^n$ is a \colb{closed and convex set with a non-empty interior}. For the sake of brevity, henceforth we denote
\begin{align}
g_t(\x) := h_t(\x) + 1_{\X_t}(\x) := \begin{cases} h_t(\x) & \x \in \X_t \\
\infty & \x \notin \X_t
\end{cases}
\end{align} 
and $f_t(\x) := \frac{1}{N}\sum_{i=1}^N f_t^i(\x)$ so that \eqref{p} may also be written as $\x_t^\star = \arg\min_\x f_t(\x) + g_t(\x)$. In the multi-agent setting considered here, the function $f_t^i$ is private to node $i$ while $g_t$ is known at all the nodes. It is remarked that the static version of \eqref{p} has found several applications in signal processing \cite{antonello2018proximal} and learning problems \cite{polson2015}. 

Algorithms to solve \eqref{p} are developed and analyzed within the rubric of online convex optimization, where the learning process is modeled as a sequential game between the agents and an adversary \cite{zinkevich2003online}. At each time $t$, agent $i$ plays an action $\x_t^i$ and in response, receives the functions $(f_t^i, g_t)$ from the adversary. Over a horizon of $T$ times, the goal of the agent is to minimize the cumulative dynamic regret given by \cite{shahrampour2018distributed}:
\begin{align}\label{regret}
\Reg:=\frac{1}{N}\sum_{t=1}^{T}\sum_{i=1}^{N}[\ell_t(\x_t^i)-\ell_t(\x^\star_t)].
\end{align}
The dynamic regret measures the loss incurred by the agent against that of a time-varying clairvoyant. Observe that the dynamic regret is different from and more stringent than the  \emph{static} regret where the benchmark is not time-varying \cite{zinkevich2003online}. The definition in \eqref{regret} also includes the modification due to the distributed nature of the problem, first suggested in \cite{shahrampour2018distributed}. For instance, the right-hand side of \eqref{regret} includes terms of the form $f^j_t(\x^i_t)$ for $i\neq j$, and therefore, regret minimization requires the agents to collaborate.

An algorithm is said to be no-regret if $\Reg$ grows sublinearly with $T$. It is well-known that the dynamic regret cannot be sublinear unless certain regularity conditions are imposed on the temporal variations of $\x_t^\star$ \cite{besbes2015non}. In the present case, the regret bounds will be developed in terms of the path length, defined as 
\begin{align}
C_T:=\sum_{t=2}^{T}\norm{\x_{t}^\star-\x_{t-1}^\star}
\end{align}
and assumed to be sublinear. For general time-varying problems with gradient feedback, the dynamic regret of any algorithm is at least $\mathcal{O}(1+C_T)$ \cite{yang2016tracking}. Dynamic regret bounds of related algorithms capable of handling non-smooth functions are summarized in Table \ref{table}. Closely related to the proposed work, a distributed online mirror descent algorithm is discussed in \cite{shahrampour2018distributed} that is more general and can handle non-differentiable loss functions. It is remarked that the centralized algorithm proposed in \cite{zhang2016improved} is not included in Table \ref{table} since it requires the loss function to be self-concordant. 

The present work develops an online and distributed algorithm for solving \eqref{p}, where the agents minimize the regret collaboratively and without a central controller or fusion center. While the static variant of the problem can be readily solved via consensus, such an algorithm is not directly applicable to the time-varying problem at hand. Specifically, while the spread of information is limited due to the time-varying communication graph $\G_t$, the problem parameters $(f_t, g_t)$ continue to change with $t$, regardless. Towards this end, we adopt the multi-step consensus idea from \cite{chen2012fast}. Sublinear regret will be obtained by carefully balancing the additional accuracy obtained from running the consensus step for multiple times against the excess regret accumulated in the meanwhile.

\section{Proposed distributed algorithm}\label{seciii}
\subsection{Motivation}
In order to motivate the proposed algorithm, observe that since $\ell_t$ is strongly convex, the optimality condition for \eqref{p} may be written as
\begin{align}
\x_t^\star = \prox_{g_t}^\alpha(\x_t^\star - \alpha \nabla f_t(\x^\star)) \label{opt}
\end{align}
where recall that $f_t(\x^\star) = \frac{1}{N}\sum_{i=1}^N f_t^i(\x^\star)$ and the proximal operator $\prox_{g_t}^\alpha$ is defined as:
\begin{align}\label{prox}
\prox_{g_t}^{\alpha}(\x) &= \arg\min_{\u} \    g_t(\u) + \frac{1}{2\alpha}\norm{\u-\x}^2  \\
&= \arg\min_{\u \in \X_t} h_t(\u) + \frac{1}{2\alpha}\norm{\u-\x}^2
\end{align}
and $\alpha > 0$ is the step-size. In a  centralized setting, the form of the optimality condition in \eqref{opt} is suggestive of the proximal OGD algorithm \cite{dixit2018online} that takes the form 
\begin{align}
\x_{t+1} = \prox_{g_t}^\alpha(\x_t - \alpha \nabla f_t(\x_t)) \label{ipogd}
\end{align}
for all $t\geq 1$. Though the algorithm in \eqref{ipogd} is provably no-regret, it is not usable in distributed settings where the timely evaluation of the average gradient $\nabla f_t(\x_t) = \tfrac{1}{N}\sum_{i=1}^N \nabla f_t^i(\x_t)$ is challenging. Indeed, due to the dynamic and arbitrary nature of the communication graphs $\G_t$, even collecting the gradients $\nabla f_t^i(\x_t)$ from each node is not straightforward.

The consensus algorithm, where the information transmission occurs only over the edges of a given graph $\G_t$, has been widely used for distributed averaging \cite{nedic2015distributed}. While the consensus algorithm converges only asymptotically, it is still possible for the nodes to obtain an approximate version of the average in a few iterations. The idea of running the consensus for a few iterations in order to approximately calculate the average gradient $\tfrac{1}{N}\sum_{i=1}^N \nabla f_t^i(\x_t^i)$ was first proposed in \cite{chen2012fast} for static problems, and will also be utilized here. A key complication that occurs in the dynamic setting is that while agents carry out the averaging via consensus routine, the problem parameters $(f_t, g_t)$ continue to change. Therefore it becomes important to carefully plan the number of consensus steps that will result in a sufficiently high update accuracy without losing track of $\x_t^\star$.

\subsection{The Distributed Proximal OGD Algorithm} 
Building upon the classical distributed proximal-gradient algorithm \cite{chen2012fast}, the proposed distributed proximal OGD (DP-OGD) algorithm takes up multiple time slots to complete each iteration. The time-varying functions $(f_t^i,g_t)$ are sampled at times $\T:= (t_1, t_2, \ldots, t_K)$ where 
\begin{align}
t_{k+1} = t_k + S(k) + 2 \label{tdiff}
\end{align}
so the $k$-th iteration starting at time $t_k$ takes up $S(k)+2$ time slots \colb{where $S(k)$ denotes the  number of consensus steps at iteration $k$. The two additional time slots are reserved for the gradient update and for the application of the proximal operator. } In order to establish the regret bounds, we will generally choose $S(k)$ to be a non-decreasing function of $k$. Associate time-varying weights $A^{ij}_t$ with each edge $(i,j)\in\E_t$, while let $A^{ij}_t = 0$ for all $(i,j)\notin\E_t$. Let the matrix $\A_t \in \R^{N\times N}$ collect all the edge weights $\{A^{ij}_t\}_{i,j}$. The full DP-OGD algorithm starts with an initial $\{\x_1^i\}$ and consists of the following updates 
\begin{subequations}\label{tup}
\begin{align}
\z_{t+1}^i &= \x_t^i - \alpha \nabla f_t^i(\x_t^i) & t \in \T\label{ztup}\\
\z_{t+1}^i &= \sum_{j: (i,j)\in \E_t} A_t^{ij}\z_t^j	& t, t+1 \notin \T \label{ztcon}\\
\x_{t+1}^i &= \prox_{g_{\t}}^\alpha(\z_t^i) & t+1 \in \T \label{xtup}
\end{align}
\end{subequations}
where we let $\t:=\max_{\tau} \{\tau \in\T | \tau \leq t\}$. Note that only one of the three updates steps in \eqref{ztup}-\eqref{xtup} is carried out depending on the value of $t$. \colb{These updates in time index $t$ have been effectively summarized in Algorithm \ref{meta_algo}}.
 
\colb{\begin{algorithm}
	\caption{DP-OGD Algorithm in time index $t$}\label{meta_algo}
\colb{\begin{algorithmic}[1]	
		\STATE {\textbf{Initialize}} $\{\x_{1}^{i}\}_{i\in\N}$, $\{\z_{1}^{i}\}_{i\in\N}$, \textbf{step-size} $\alpha$, and $\T:= (t_1, t_2, \ldots, t_K)$ where $t_{k+1} = t_k + S(k) + 2 $
		\STATE {\textbf{for} $t=1,\cdots,T$} \textbf{do} 
		\STATE \ \ \ \ \ {\textbf{for} $i=1,2,\cdots,N$} \textbf{do} 
		\IF{$t \in \T$}
		\STATE \ \ \ \ \ \ \ \ \textbf{Get} exact gradient $\nabla {f}_t^{i}(\x_t^{i}) $ and function  ${g}_t$	  		
		\STATE \ \ \ \ \ \ \ \ \textbf{Update} $ {\z}_{t+1}^{i} = \x_t^{i} - \alpha \nabla {f}_t^{i}(\x_t^{i}) $ 
		\ENDIF
		\IF{$t, t+1 \notin \T$}
		\STATE \ \ \ \ \ \ \ \ \textbf{Receive} $\z_t^j$ from all $j: (i,j)\in \E_t$  
		\STATE \ \ \ \ \ \ \ \ \textbf{Update} $  \z_{t+1}^i = \sum_{j: (i,j)\in \E_t} A_t^{ij}\z_t^j $  
		\ENDIF
		\IF{$t+1 \in \T$}
		\STATE \ \ \ \ \ \ \ \ \textbf{Update} $  \x_{t+1}^i = \prox_{g_{\t}}^\alpha(\z_t^i)  $ \STATE \ \ \ \ \ \ \ \ \ \ where $ \t:=\max_{\tau} \{\tau \in\T | \tau \leq t\}$
		\ENDIF
		\STATE \ \ \ \ \ \textbf{end for} 
		\STATE \textbf{end for}
	\end{algorithmic}}
   \end{algorithm}}	

\begin{figure*}
	\centering
	\vspace{-5mm}
	\includegraphics[scale=0.5]{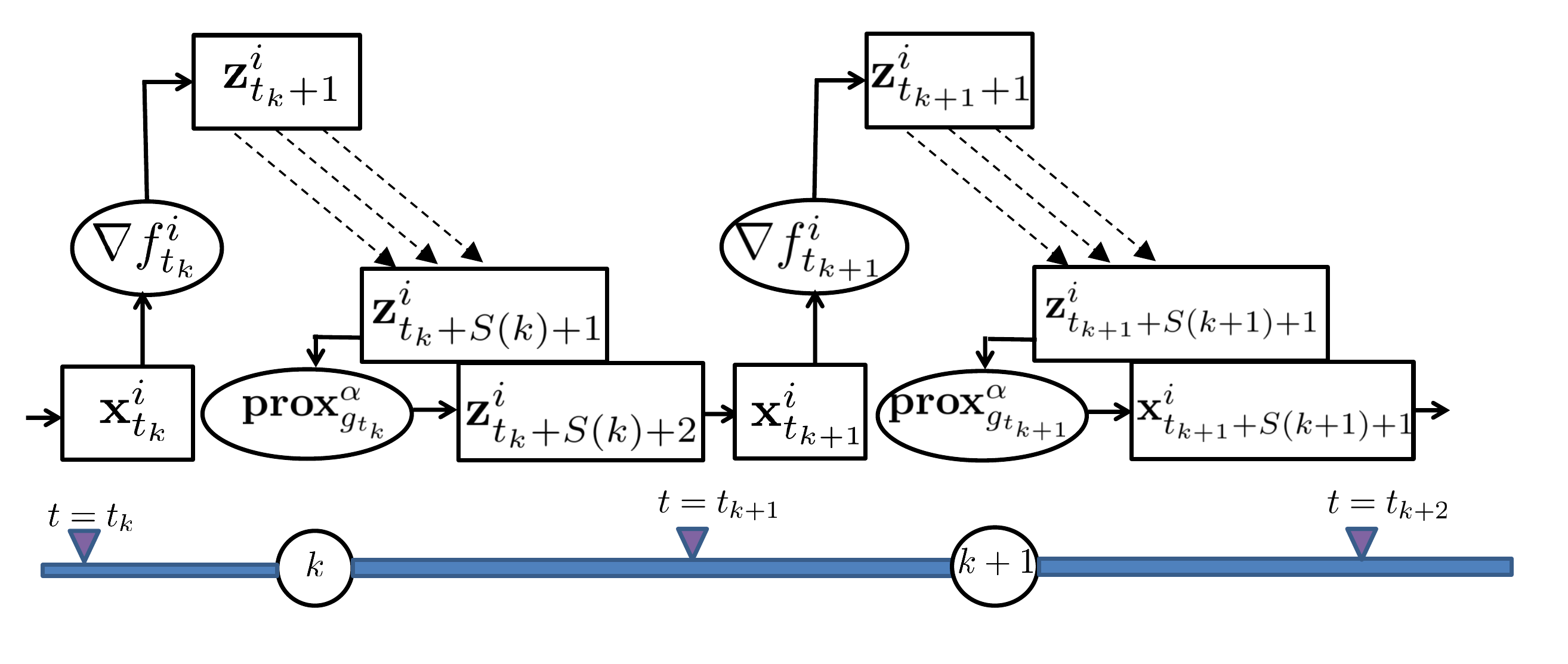}
	\caption{\text{DP-OGD updates at node $i$.}}
	\label{fig:fig1_schema}
\end{figure*}

In order to better understand the proposed algorithm, it may be instructive to write down the updates in \eqref{tup} in terms of the iteration index $k$:
\begin{subequations}
\begin{align}
\z_{t_k+1}^i &= \x_{t_k}^i - \alpha \nabla f_{t_k}^i(\x_{t_k}^i) \label{ztkup}\\
\z_{{t_k}+s+1}^i &= \sum_{j: (i,j)\in \E_{{t_k}+s}} A_{{t_k}+s}^{ij}\z_{{t_k}+s}^j \hspace{5mm} 1 \leq s\leq S(k) \label{ztkcon}\\
\x_{t_{k+1}} &= \x_{{t_k}+S(k)+2}^i = \prox_{g_{t_k}}^\alpha(\z_{{t_k}+S(k)+1}^i) \label{xtkup}
\end{align}
\end{subequations}
For the sake of brevity, we introduce new (capped) variables and functions whose subscript indicates the iteration index instead of the time index. Specifically, for all $k=1, \ldots, K$, let
\begin{subequations}\label{konly}
\begin{align}
\xh^i_k &:= \x^i_{t_k} & \zh^i_k &:= \z^i_{t_k+1} \\
\fh_k^i &:= f_{t_k}^i & \gh_k &:= g_{t_k} \\
\yh^i_k &:= \z^i_{t_k+S(k)+1} 
\end{align}
\end{subequations}
Likewise, the optimum at time $t_k$ is specified as $\xh^\star_k$. With the iteration-indexed notation in \eqref{konly}, the DP-OGD updates may simply be written as
\begin{subequations}
\begin{align}
\zh_k^i &= \xh_k^i - \alpha \nabla \fh_k^i(\xh_k^i) \label{zkup}\\
\yh_k^i &= \sum\nolimits_{j=1}^N W_k^{ij}\zh_k^j	 \label{zkcon}\\
\xh_{k+1}^i &= \prox_{\gh_k}^\alpha(\yh_k^i) \label{xkup}
\end{align}
\end{subequations}
where $\colb{Q}_k^{ij}$ denotes the $(i,j)$-th entry of the  matrix 
\begin{align}
\colb{\Q}_k &:= \A_{t_k+S(k)}\ldots\A_{t_k+1}.
\end{align}
Here, the $S(k)$-steps of consensus are all condensed into a single equation \eqref{zkcon}. In order to better understand \eqref{zkcon}, collect the variables $\{\zh_k^i\}$, $\z_k^i$, and $\yh_k^i$ into super-vectors $\zh_k$, $\z_k$, and $\yh_k$ respectively. Recursively applying the consensus averaging, it follows that 
\begin{align}
\yh_k &= \z_{{t_k}+S(k)+1}= \A_{t_k+S(k)}\z_{t_k+S(k)} \nonumber\\
 \ldots &= \A_{t_k+S(k)}\ldots\A_{t_k+1}\z_{t_k+1} = \colb{\Q}_k\zh_k 
\end{align}
The full algorithm is summarized in Algorithm \ref{algo_1} and the updates at node $i$ are shown in Fig. \ref{fig:fig1_schema}. It is remarked that Algorithm \ref{algo_1} is still conceptional since we have left $S(k)$ unspecified. An appropriate choice of $S(k)$ is necessary to obtain a tight regret bound and a detailed discussion for the same will be provided in Sec. \ref{perfa}. 

	  \begin{algorithm}
	  	\caption{Distributed Proximal Online Gradient Descent}\label{algo_1}
	  	\begin{algorithmic}[1]
	  		\STATE {\textbf{Initialize}} $\{\hat\x_{0}^{i}\}_{i\in\N}$	and step-size $\alpha$ 
	  		\STATE {\textbf{for} $k=0,1,\cdots,K$} \textbf{do} 
	  		\STATE \ \ \ \ \ {\textbf{for} $i=1,2,\cdots,N$} \textbf{do} 
	  		\STATE \ \ \ \ \ \ \ \ \textbf{Get} exact gradient $\nabla \hat{f}_k^{i}(\hat\x_k^{i}) $ and function  $\hat{g}_k$	  		
	  		\STATE \ \ \ \ \ \ \ \ \textbf{Update} $ \hat{\z}_k^{i} = \hat\x_k^{i} - \alpha \nabla \hat{f}_k^{i}(\hat\x_k^{i}) $ 
	  		 \STATE \ \ \ \ \ \textbf{end for} 
	  		\STATE \ \ \ \ \ {\textbf{for} $i=1,2,\cdots,N$} \textbf{do} 
	  		\STATE \ \ \ \ \ \ \ \ \textbf{Update} $  \yh_k^{i} = \sum_{j=1}^{N} \colb{Q}_k^{ij}\zh_k^{j} $      
	  	    \STATE \ \ \ \ \ \ \ \ \textbf{Update} $  \xh_{k+1}^{i} = \prox_{\gh_k}^{\alpha}(\yh_k^{i})  $
	  		\STATE \ \ \ \ \ \textbf{end for} 
	  			  		\STATE \textbf{end for }
	  	\end{algorithmic}
	  		  	\label{algo1}
	  \end{algorithm}

\section{Performance Analysis}\label{perfa}
This section develops the regret rate for the proposed algorithm. As already discussed, each iteration involves the gradient update at time $t_k$, $S(k)$ consensus steps, and a proximal update at time $t_{k+1}-1$. As compared to the existing proximal OGD algorithm, the analysis of DP-OGD algorithm is complicated due to the additional regret that is incurred from subsampling the functions $(f_t, g_t)$ at times $t\in\T$. We begin with discussing some preliminaries before proceeding to the assumptions and the regret bounds.

\subsection{Preliminaries}
The regret rate of Algorithm \ref{algo_1} will be analyzed using the network averages at each iteration $k$, defined as
\begin{align}\label{netav}
\xb_k =& \frac{1}{N}\sum_{i=1}^{N}\xh_k^{i} & \zb_k =& \frac{1}{N}\sum_{i=1}^{N}\zh_k^{i}.
\end{align}

Towards establishing the regret rates, we begin with casting the proposed algorithm as an inexact proximal gradient algorithm that can be viewed as a generalized version of \cite{dixit2018online}. Using the network averages defined in \eqref{netav}, it is possible to write \eqref{zkup} as an inexact gradient update step:
\begin{align}
\zb_k &= \xb_k - \frac{\alpha}{N}\sum_{i=1}^N \colb{\nabla \fh_k^i(\xh_k^i)} \label{izup0} \\
&= \xb_k - \alpha[\nabla \fh_k(\xb_k) + \e_k] \label{izup}
\end{align}
where recall that $\fh_k(\x) := \frac{1}{N}\sum_{i=1}^N \fh_k^i(\x)$ and 
\begin{align}\label{error_Definition}
\e_k =  \frac{1}{N}\sum_{i=1}^N(\nabla \fh_k^i(\xh_k^i) - \nabla \fh_k^i(\xb_k) ).
\end{align}
Along similar lines, we write the inexact proximal update as
\begin{align}
\xb_{k+1} &= \prox_{\gh_k,\epsilon_k}^{\alpha}(\zb_k) \\
&=:\prox_{\gh_k}^{\alpha}(\zb_k) + \ve_k. \label{ixup}
\end{align}
\colb{where the $\epsilon_k$-proximal operator $\prox_{\gh_k,\epsilon_k}^{\alpha}$ is as defined in \cite{rates11schmidt} which implies that
\begin{align}
\frac{1}{2\alpha}&\norm{\xb_{k+1}-\zb_k}^2 + \gh_k(\xb_{k+1}) \nonumber\\ 
&\leq \epsilon_k + \min_{\x \in \mathcal{\X}}\bigg(\frac{1}{2\alpha}\norm{\x-\zb_k}^2 + \gh_k(\x)\bigg).
\end{align}} As compared to \eqref{zkcon}-\eqref{xkup}, the error incurred in using the inexact proximal operation can be expressed as
\begin{align}\label{ve}
\ve_k = \frac{1}{N}\sum_{i=1}^N\prox_{\gh_k}^\alpha(\y_k^i)-\prox_{\gh_k}^{\alpha}(\zb_k).
\end{align}
for all $k\geq 1$. At this stage, \colb{$\ve_k$} is left unspecified and appropriate bounds will be developed later. Having written the proposed DP-OGD algorithm as an inexact proximal OGD variant, the rest of the analysis proceeds as follows: (a) development of bounds for updates in \eqref{izup}-\eqref{ixup} in terms of $\norm{\e_k}$ and $\norm{\ve_k}$; (b) development of bounds on $\norm{\e_k}$ and $\norm{\ve_k}$; and finally (c) substitution of these bounds to obtain the required regret bounds. It is remarked that such an approach is flexible and readily extendible to other variants where the sources of gradient and proximal errors may be different; e.g., quantization, asynchrony, or computational errors.

\subsection{Assumptions}
The assumptions required for developing the regret bounds are discussed subsequently. The first three assumptions pertain to the properties of functions $f_t^i$ and $g_t$. 

\begin{assum}\label{sm}
The functions $f_t^i$ are $L$-smooth, i.e., for all $\colb{\x,\y\in\X_t}$, it holds that
\begin{align}\label{asmooth}
\norm{\nabla f_t^{i}(\x) - \nabla f_t^{i}(\y)}\leq L\norm{\x-\y}  
\end{align}
for all $t\geq 1$ and $1\leq i\leq N$.  
\end{assum}

\begin{assum}\label{lip}
The cost functions $f_t^i$ and $g_t$ are Lipschitz continuous, i.e., for all \colb{$\x, \y \in \X_t$}, we have that
\begin{align}
\norm{f_t^i(\x)-f_t^i(\y)} &\leq M\norm{\x-\y} \label{lipf}\\
\norm{g_t(\x)-g_t(\y)} &\leq M\norm{\x-\y} \label{lipg}
\end{align}
where we have used the same Lipschitz constant for $\{\{f_t^i\}_{i=1}^N, g_t\}_{t=1}^T$ for the sake of brevity and without loss of generality. Note that \eqref{lipf}-\eqref{lipg} also imply that the corresponding (sub-)gradients are bounded by $M$. 
\end{assum}

\begin{assum}\label{con}
The functions $f_t^i$ are $\mu$-convex, i.e., for all $\x, \y \in \X_t$, it holds that
\begin{align}
\ip{\nabla f_t^i(\x)-\nabla f_t^i(\y),\x-\y} \geq \mu\norm{\x-\y}^2 
\end{align}
\end{assum}

It is remarked that Assumptions \ref{sm}-\ref{con} are standard and applicable to a wide range of problems arising in machine learning, signal processing, and communications. The present analysis depends critically on these assumptions and the required regret bounds only hold when $\mu >0$ and the parameters $L$ and $M$ are bounded. 
The next two assumptions pertain to the network connectivity and the choice of weights.  
\begin{assum}\label{bcon}
The graph $\G_t$ is $B$-connected, i.e., there exists some $B \in \mathbb{N}$ such that the graph 
\begin{align}\label{graph}
\G_t^B:=\left(\N, \bigcup\limits_{\tau =t B}^{(t+1)B-1} \E_\tau\right)
\end{align}
is connected for all $t\geq 1$. 
\end{assum}

\begin{assum}\label{wt}
The weight matrices $\{\A_t\}_{t=1}^T$ have positive entries and satisfy the following properties for each $t \geq 1$:
\begin{enumerate}
\item \emph{Double stochasticity:} it holds that $\sum_{i=1}^N A_t^{ij} = \sum_{j=1}^N A_t^{ij} = 1$;
\item \emph{Lower bounded non-zero entries:} there exists $\eta > 0$ such that $A_t^{ij} \geq \eta$ for all $(i,j) \in \E_t$;
\item \emph{Non-zero diagonal entries:} it holds that $A_t^{ii} > \eta$ for all $i \in \N$.
\end{enumerate}
\end{assum}

Assumptions \ref{bcon}-\ref{wt} imply that the entries of the weight matrix $\colb{\Q}_k = \A_{t_k+S(k)}\ldots\A_{t_k+1}$ are close to $1/N$ in the following sense \cite[Proposition 1(b)]{nedic2009distributed}:
\begin{align}\label{wbound}
\abs{\colb{Q}_k^{ij}-\frac{1}{N}} \leq \Gamma \gamma^{S(k)-1}
\end{align}
where defining $\omega := \eta^{(N-1)B}$, we generally have that $\Gamma = 2\frac{\omega+1}{\omega(1-\omega)}$ and $\gamma = (1-\omega)^{\tfrac{1}{B}}$. Inequality \eqref{wbound} holds the key to obtaining fast consensus over time-varying graphs and as we shall see later, the regret rate of the algorithm would depend critically on the value of $\gamma$. 

For the next assumption, let the total number of time slots required to carry out $K$ updates be given by
\begin{align}
\hat{S}_K:=\sum_{k=1}^K [S(k)+2] = T \label{hats}
\end{align}
where recall that $S(k)$ is the number of consensus steps at the $k$-th iteration. Inverting the relationship \eqref{hats}, for any $T$, it holds that $K = \max\{K | \hat{S}_K \leq T\} =: S_T$. The following assumption is required to ensure that the regret bounds are sublinear. 

\begin{assum}\label{sk}
The number of consensus steps $S(k)$ is a non-decreasing function of $k$ and the number of consensus steps at the last iteration $S(K) = S(S_T) =: R_T$ is sublinear in $T$.
\end{assum}

\colb{An implication of Assumption \ref{sk} is that there cannot be too many consensus steps at any iteration $k$. } Having stated all the required assumptions, we are now ready to state the main results of the paper. 

\subsection{Regret Bounds}\label{secreg}
As already discussed, we begin with developing some bounds for the generic inexact proximal gradient method \eqref{izup}-\eqref{ixup}. Bounds on the error incurred from using the proposed distributed algorithm will be developed next. The final regret bounds would follow from combining these two results. The section concludes with some discussion on the nature of the bounds for various choices of $S(k)$. 

The first lemma bounds the per-iteration progress of the iterate $\xb_k$ in terms of its change in distance from the current optimum $\xh_k^\star$. 
\begin{lem}\label{lem1}
Under Assumptions \ref{sm}-\ref{con}, the updates in \eqref{izup}-\eqref{ixup} satisfy
\begin{align}\label{lemi}
&\norm{\xb_{k+1} - \xh_k^\star} \leq \rho \norm{\xb_k - \xh_k^\star} + \delta_k
\end{align}
where, $\rho^2:= 1 + \alpha^2 L^2 - 2 \alpha\mu$ and $\delta_k:= \norm{\ve_k} + \alpha\norm{\e_k}$.
\end{lem}
The proof of Lemma \ref{lem1} is provided in Appendix \ref{proof_first} and utilizes the strong convexity and smoothness properties of $\fh_t$ as well as the triangle inequality. The result in Lemma  \ref{lem1} states that the distance between $\xb_{k+1}$ and the current  optimal $\xh_k^\star$ is less than a $\rho$-fraction of the distance between $\xb_k$ and $\xh_k^\star$, but for an error term $\delta_k$ that arises due to $\e_k$ and $\ve_k$ in the updates steps. It is remarked that the result in \eqref{lemi} subsumes all existing results in \cite{bedi2018tracking,dixit2018online}. 

Taking summation over $k\geq 1$ and rearranging, we obtain the following corollary, whose proof is deferred to Appendix \ref{proof_first}. 
\begin{cor}  \label{cor1}	
Under Assumptions \ref{sm}-\ref{con} and for $0 < \alpha < 2\mu/L^2$, the updates in \eqref{izup}-\eqref{ixup}  satisfy 
\begin{align}\label{first_bound}
\sum_{k=1}^K\norm{\xb_k- \xh_k^\star}  \leq & \frac{\rho}{1-\rho}\norm{\xb_{0} - \xh_{0}^\star} + 	\sum_{k=1}^K\frac{\norm{\xh_k^\star - \xh_{k-1}^\star}}{1-\rho}\nonumber\\& 	  +  \frac{1}{1-\rho}\sum_{k=0}^K \delta_k.
\end{align}	  		 	
\end{cor}	  		 
 The result in Corollary \ref{cor1} provides an upper bound on the  cumulative  distance between the average current iterate $\xb_k$ and optimal $\xh_k^\star$ over $K$ time instances. \colb{Next, Lemma \ref{cor2} develops a simple bound on the iterate norm with the proof provided in Appendix \ref{proofcor2}}. 
 
\begin{lem}\label{cor2}	
Under Assumption \ref{lip}, \ref{bcon}, and \ref{wt}, the iterates generated by the DP-OGD algorithm are bounded as
\begin{align}
	\sum_{i=1}^{N}\norm{\xb_k- \xh_k^{i} } &\leq  2  \Gamma \gamma^{S(k-1)-1} N \sum_{i=1}^{N}\norm{\zh_{k-1}^i }\label{cor2a}\\ 
 \sum_{i=1}^{N}\norm{ \zh_{k-1}^{i}} &\leq  \sum_{i=1}^{N}\norm{\zh_{0}^{i}} +  2\alpha N  M (k-1)   \label{cor2b}   
 \end{align}	
\end{lem}	
The bounds in Lemma \ref{cor2} are not surprising given that the size of each update step is bounded implying that after $k$ steps, none of the iterates can be more than $\mathcal{O}(k)$ far from the starting point. Recall that the inexact proximal OGD algorithm updates in \eqref{izup}-\eqref{ixup} are really the DP-OGD updates in disguise, with specific definitions of $\e_k$ and $\ve_k$. The next lemma provides a convenient bound on the error term $\delta_k$ that will subsequently be used to obtain the regret bounds. 
  		 
\begin{lem}\label{lemdel}
Under Assumptions \ref{sm}-\ref{wt}, the error $\delta_k=\norm{\ve_k}+\alpha\norm{\e_k}$ is bounded as
  	\begin{align}  
\delta_k &\leq  2\alpha L \Gamma \gamma^{S(k-1)-1} \left(\sum_{i=1}^{N}\norm{\zh_{0}^{i}} +  2\alpha N  M  (k-1)\right)   \nonumber \\
  		 & +\Gamma \gamma^{S(k)-1} \left(\sum_{i=1}^{N}\norm{\zh_{0}^{i}} +  2\alpha N  M  k    \right)   \label{delta_k}
 \end{align}
\end{lem}
 The proof of Lemma \ref{lemdel} is provided in Appendix \ref{lemma_2_proof}. It provides a bound on the error sequence generated by the inexact algorithm in terms of a geometric-polynomial sequence. The first term of \eqref{delta_k} results from the gradient error $\e_k$ while the second term bounds the norm of the proximal error $\ve_k$. 

We are now ready to convert the results obtained in terms of the iteration index $k$ into bounds that depend on time index $t$. Recall that since the $k$-th iteration incurs $S(k)+2$ time slots, the last iteration incurs $R_T := S(K)$ time slots where $K := S_T$, as stated in Assumption \ref{sk}. Also let $E_T := \sum_{k=1}^{S_T} \gamma^{S(k)}k$ so that the initial bounds can be developed in terms of $R_T$ and $E_T$. Specific examples of $S(k)$ will subsequently be provided to yield bounds as explicit functions of $T$. As a precursor, consider a simple example when $S(k) = k$, so that $R_T = K = \mathcal{O}(\sqrt{T})$ and $E_T = \mathcal{O}(1)$. Next, the following lemma reconciles the two definitions of the path length.
\begin{lem}\label{pl}
It holds that 
\begin{align}
\sum_{k=1}^K\norm{\xh_k^\star - \xh_{k-1}^\star} &\leq  C_T   \label{iteration_ind_2}
\end{align}
where $K$ is such that $\sum_{k=1}^K (S(k)+2) = T$. 
\end{lem}
\begin{IEEEproof} The result follows from the use of triangle inequality:
	\begin{align}
&\sum_{k=1}^K\norm{\xh_k^\star - \xh_{k-1}^\star} =  \sum_{k=1}^K\norm{\x_{t_k}^\star - \x_{t_{k-1}}^\star} \nonumber \\ 
&\leq\sum_{k=1}^K\sum_{\tau = t_{k-1}}^{t_k-1} \norm{\x^\star_{\tau+1}-\x^\star_\tau}= \sum_{t=1}^{T}\norm{\x_{t}^\star - \x_{t-1}^\star} = C_T \nonumber.
	\end{align}
\end{IEEEproof}

Note that in order to calculate the dynamic regret in \eqref{regret}, it is necessary to define $\x_t$ for all $t$. Towards this end, let $\x_t^i = \x_{\t} = \x_{t_k}^i$ for $k$ such that $t_k \leq t < t_{k+1}$. Finally, we provide the main result of the paper in the following Theorem.

\begin{thm}\label{thm1}
Under Assumptions \ref{sm}-\ref{sk} and $\tfrac{1}{T} \ll \alpha < \frac{2\mu}{L^2}$, the proposed DP-OGD algorithm incurs the following dynamic regret
\begin{align}\label{rbound}
\Reg\leq \mathcal{O}(R_T(1+E_T+C_T)).
\end{align}  
\end{thm}
The proof of Theorem \ref{thm1} is provided in Appendix \ref{Proof_thm_1}.
Here, the regret bound is worse than $C_T$ since the algorithm updates are sporadic resulting in additional factor of $R_T$. The result in \eqref{rbound} depends on the choice of $\{S(k)\}$ through $E_T$ and $R_T$. We now discuss the explicit form of the regret bound for a few choices of $S(k)$. 

\subsubsection{Logarithmically increasing $S(k)$} Consider the case 
\begin{align}
S(k) = \lfloor c\log(k) \rfloor
\end{align}
where $c > 1$ is left unspecified at this stage. Given $T$, the number of iterations are given by the largest $K$ that satisfies the inequality 
\begin{align}
\sum_{k=1}^K (\lfloor c\log(k) \rfloor + 2) \leq T
\end{align}
For the sake of brevity, let $T \gg 1$ and likewise $K \gg 1$ so that only the dominant terms may be retained. Ignoring the floor function and using Stirling's approximation \cite{whittaker1996course}, we have that
\begin{align}
T = \mathcal{O} (c K \log K)
\end{align}
or equivalently, 
\begin{align}\label{new_ref}
K = S_T = {\mathcal{O}\left(\exp\left(W\left(\tfrac{T}{c}\right)\right)\right)}
\end{align}
where {$W$} denotes the Lambert {$W$} function. It follows that $R_T = S(S_T) = \mathcal{O}({cW(\tfrac{T}{c})})$. Also note that 
\begin{align}
E_T &= \sum_{k=1}^{S_T} \gamma^{\lfloor c\log(k)\rfloor + 2} k \approx \sum_{k=1}^{S_T} \gamma^2 k^{c\log (\gamma) + 1 } \\
&=\mathcal{O}(S_T^{c\log(\gamma)+2}) \\
& \approx \mathcal{O}(T^{c\log(\gamma)+2})
\end{align}
where the last step uses the approximation {$W(x) \approx \log(x)$} for large $x$. The overall regret bound thus becomes
\begin{align}
\Reg \leq \mathcal{O}(c\log(T)(1+T^{c\log(\gamma)+2}+C_T))
\end{align}
For the regret to be sublinear, it is necessary but not sufficient that $C_T$ is also sublinear. For instance, if $C_T = T^\beta$ with $\beta<1$, it is also required to hold that $c\log(\gamma)+2 < \beta$ or equivalently one must choose $c > -(2-\beta)/\log(\gamma)$ so that the $C_T$ term dominates the summation. In this case, $\gamma$ is not allowed to be arbitrarily close to 1. Instead, it is required that $c \log(T)< T^{1-\beta}$ or equivalently, we have $\log(\tfrac{1}{\gamma}) > \tfrac{(2-\beta)\log{(T)}}{T^{1-\beta}}$. In other words, for a given $\gamma$, it is always possible to choose an appropriate value of the parameter $c$, though the regret bound will only be $\mathcal{O}(\log(T)(1+C_T))$ for sufficiently large $T$. 

\subsubsection{Constant $S(k)$} Taking $S(k) = \lfloor T^u \rfloor $ for all $k \geq 1$, it is required that
\begin{align}
&\sum_{k=1}^K (\lfloor T^u  \rfloor + 2) = T \\
\Rightarrow & S_T = K \approx \frac{T}{T^u+2} \approx T^{1-u}
\end{align}
for $K$ and $T$ sufficiently large. In this case,  $R_T = T^u$ and 
\begin{align}
E_T  = \sum_{k=1}^{T^{1-u}} \gamma^{T^u} k \\
\approx \gamma^{T^u}T^{2-2u}
\end{align}
yielding the final regret bound
\begin{align}
\Reg \leq \mathcal{O}(T^u(1+\gamma^{T^u}T^{2-2u}+C_T))\label{regc}
\end{align}
In order to write the regret in explicit form, let $C_T = T^\beta$ for some $\beta \in [0,1)$. Then the regret in \eqref{regc} is sublinear when $u+\beta < 1$ or $u < 1-\beta$. However $u$ cannot be arbitrarily small or else the term $\gamma^{T^u}T^{2-2u}$ would become too large. The minimum value of the regret is obtained for the choice of $u$ such that 
\begin{align}
T^\beta &= \gamma^{T^u}T^{2-2u}   \\
T^u\log(\tfrac{1}{\gamma}) &= (2-2u-\beta)\log(T) \\
\Rightarrow \hspace{1cm} u &= \frac{(2-\beta)}{2} - \frac{{W}\bigg(T^{\frac{(2-\beta)}{2}}\log{(1/\gamma)^{1/2}}\bigg)}{\log{(T)}}
\end{align}
where we have used the result from  \cite{corless1996lambertw}. Since $T$ is large, we make use of the approximation ${W(x)} \approx \log x - \log\log x$, which yields
\begin{align}
u  &\approx \frac{\log{\bigg(\frac{2}{\log{(1/\gamma)}}\bigg)}}{\log{(T)}} + \frac{\log{\bigg(\log{\bigg(T^{\frac{(2-\beta)}{2}}\log{(1/\gamma)^{1/2}}\bigg)}\bigg)}}{\log{(T)}} \nonumber\\
T^u &\approx  \frac{1}{\log{(\tfrac{1}{\gamma})}} \left(2\log{\log{(1/\gamma)^{1/2}}} + (2-\beta)\log T\right) \nonumber
\end{align}
Therefore the regret bound can be approximately written as
\begin{align}
\Reg &= \mathcal{O}\bigg( T^u(1+C_T)\bigg)  \approx \mathcal{O}\bigg(\log{(T)}(1+C_T)\bigg)\label{regc2}
\end{align} 
which is sublinear as long as $\gamma$ is not too close to 1 and $\beta<1$. As in the previous case, the optimal choice of $u$ still requires a sufficiently large value of $T$. \colb{In summary, the step-size may be chosen as $\tfrac{1}{T} \ll \alpha < \frac{2\mu}{L^2}$ while $S(k)$ may be chosen as either $\lfloor c\log(k) \rfloor$ or $\lfloor T^u \rfloor$ in Algorithm \ref{algo1}. }

Recall that the centralized proximal OGD algorithm achieves a regret of $\mathcal{O}(1+C_T)$ which is also optimal for any online algorithm; see \cite{yang2016tracking}. Remarkably, the dynamic regret of the DP-OGD is only worse by a $\log(T)$ factor, possibly arising out of the distributed operation over an intermittently connected graph. 


\section{Numerical results}\label{numerical}
The performance of the proposed DP-OGD algorithm is tested for the dynamic sparse signal recovery problem where the goal is to estimate a time-varying sparse parameter.  Such problems have been widely studied in literature and can be broadly classified into those advocating adaptive filtering-based algorithms, those formulating the problem within the sparse Bayesian learning framework, and finally those utilizing tools from dynamic or online convex optimization.

\subsection{The Dynamic Sparse Recovery Problem}
Consider a WSN with $N$ sensors connected over a time-varying graph $\G_t$. The parameter of interest is a time-varying sparse signal $\u_t \in \mathbb{R}^n$. At time $t$, sensor $i$ makes $d$ measurements according to the following model
\begin{align}
\textbf{y}^i_t = \C^i_t\u_t + \textbf{v}^i_t 
\end{align}
where $\C^i_t \in \R^{d \times n}$ is the observation matrix and $\v_t^i \in \R^{d}$ is the noise with unknown statistics. The online learning model detailed in Sec. \ref{Prob_for} is considered and the quantities $\{\y^i_t, \C_t^i\}_{i\in\N}$ are revealed in a sequential manner. Observe that given no other information, tracking the original parameter $\u_t$ is impossible. Instead, we settle for tracking the Elastic Net estimator of $\u_t$ given by 
 \begin{align}\label{en}
\x^{\star}_t = \argmin_{\x} \frac{1}{N}\sum_{i=1}^{N}\norm{\y_t^i-\C_t^i\x}_2^2 + \lambda \!\norm{\x}_2^2 \!+\! \sigma  \norm{\x}_1
 \end{align}
where $\lambda$ and $\sigma$ are regularization parameters. \colb{While the problem is unconstrained, gradient boundedness can be ensured by imposing a norm-ball constraint with a large radius. }


Numerical tests are carried out for a network with $N = 100$ nodes where each node makes $d = 4$ measurements in order to track a parameter of dimension $n = 50$. For the purpose of initialization, $\u_0$ is chosen to be a sparse random vector with $10$ non-zero entries. Let $\S_t:=\{n \mid [\u_t]_{n} > 0\}$ be the support of the vector $\u_t$. For time $t\geq 2$ support is updated as follows:
\begin{align}
\S_{t+1} = \begin{cases}
\S_t & \text{with probability }1-1/t \\
\{\S_t\setminus\{i_t\}\} \cup \{i'_t\} & \text{with probability } 1/t 
\end{cases}
\end{align} 
where $i_t$ is randomly chosen from the $\S_t$ and $i'_t$ is randomly chosen from the zero locations $\S^c_t$. Subsequently, noise is added to the entries of $\u_t$ and normalization is performed so that  
\begin{align}
\u_{t+1} = \frac{\u_t + \mathbf{n}_t}{\norm{\u_t + \mathbf{n}_t}}
\end{align}
where $[\mathbf{n}_t]_n \sim \mathcal{N}(0,1/t^2)$ for $n\in\S_{t+1}$ and $[\mathbf{n}_t]_n = 0$ otherwise. In other words, both the support and the non-zero values are time-varying but the variations decay over time. 

For the Elastic Net estimator, we set the parameters $\sigma = 0.01/d^2N^2$ and $\lambda = 0.05/dN$. The parameters manually selected to ensure that $\x_t^\star$ remains close to $\u_t$ for all $t$. Since the goal of the DP-OGD, proximal OGD, and ADMM algorithms is to track  $\x_t^\star$, the resulting dynamic regret is not very sensitive to the choice of the parameters $\lambda$ and $\sigma$, and similar results can be obtained for other settings as well. 

The results of the DP-OGD algorithm depend on the choice of the weight matrices $\A_t$, which can be generated in different ways. A simple choice corresponds to the complete graph with all weights set to $1/N$. More generally, we utilize Birkhoff-von Neumann theorem \cite{dufosse2018further} to generate random doubly stochastic matrices. Specifically, let $\mathbf{P}^0:=\mathbf{I}_{N \times N}$ and generate $N-1$ matrices $\{\mathbf{P}^{j}\}_{j=1}^{N-1}$ by randomly permuting its rows. In general it holds that any convex combination of $\{\mathbf{P}^j\}$ is doubly stochastic, allowing us to use the following consensus weights
\begin{align}
\A_t^{(\iota)} &= \sum_{j=0}^{\iota} \omega_{j}\mathbf{P}^j 
\end{align}
for $\iota \geq 1$. Here, the weights must be selected so as to ensure Assumption 5 with $\eta = 2/N$. For the simulations, we utilize four kinds of graphs, corresponding to $\iota = 1$, $2$, $3$, and $N-1$. When $\iota = 1$, the resulting graph is sparse and disconnected, and the consensus entails each node sending its update to only one other node. On the other hand, the graph may often be fully connected when $\iota = N-1$.

\subsection{Dynamic Sparse Recovery Algorithms} For the tests, we consider three different algorithms: the proposed DP-OGD algorithm, the centralized ADMM algorithm from \cite{simonetto2015non}, and the centralized proximal OGD algorithm \cite{dixit2018online}. 

\subsubsection{Tracking via DP-OGD} Within the online setting considered here, the problem parameters $\{\y^i_t, \C_t^i\}_{i\in\N}$ are revealed after the iterates $\{\x_t^i\}_{i\in\N}$ have been obtained. Since the objective function is strongly convex but has a non-differentiable component and the sensors are connected over a time-varying graph $\G_t$, we make use of the proposed DP-OGD algorithm to track $\x_t^\star$. Recalling the form of the original optimization problem, it follows that 
\begin{align}
f_t^i(\x) &= \norm{\y_t^i-\C_t^i\x}_2^2 + \lambda\norm{\x}_2^2 \label{fti}\\
g_t(\x) &= \sigma\norm{\x}_1 \label{gt}
\end{align}
which allows us to write down the updates in Algorithm \ref{algo1}. For the tests, we use $\alpha = 0.5$, which was the largest value  of $\alpha$ for which the algorithm did not diverge.

\subsubsection{Tracking via centralized 'slowed' ADMM} For the purposes of comparison, we consider the dynamic ADMM algorithm proposed in \cite{simonetto2015non} which also allows non-differentiable and time-varying objective functions. The algorithm in \cite{simonetto2015non} is however centralized and processes all the measurements in one shot. In order to carry out a fair comparison, we consider a 'slowed' version of the ADMM algorithm that has the same number of iterations per time-instant as the proposed DP-OGD algorithm. That is, instead of carrying out one iteration per time slot, the slowed ADMM is run with one iteration for $t_k \leq t < t_{k+1}$, the underlying assumption being that the intermediate time-slots are used for exchanging the measurements between the nodes. 

In order to derive the updates of the ADMM algorithm, we introduce a new variable $\z$ and reformulate the problem at time $t$ as
\begin{align}
\min_{\x,\z} &\frac{1}{N} \sum_{i=1}^N f_t^i(\x) + g_t(\z) \\
&\text{s. t. } \x = \z
\end{align}
where $f_t^i$ and $g_t$ are as defined in \eqref{fti}-\eqref{gt}. Associating a dual variable $\v$ with the constraint, the Lagrangian is given by
\begin{align}
L_t(\x, \z, \v) = f_t(\x) + g_t(\z) + \ip{\v, \x-\z} + \frac{\varrho}{2}\norm{\x-\z}_2^2 \nonumber
\end{align}
where $f_t(\x) := \tfrac{1}{N}\sum_{i=1}^N f_t^i(\x)$. Given the iterates $(\x_t, \z_t, \v_t)$ at time $t$, the ADMM updates are given by \cite{simonetto2015non}:
\begin{subequations}\label{admmup}
\begin{align}
\x_{t+1} = & \argmin_{\x} \mathcal{L}_t(\x,\z_t,\v_t) + \frac{\varpi}{2} \norm{\x - \x_t}^2\\
\z_{t+1} = & \argmin_{\z} \mathcal{L}_t(\x_{t+1},\z,\v_t) + \frac{\varpi}{2} \norm{\z - \z_t}^2\\
\v_{t+1} = & \v_t + \varrho (\x_{t+1}- \z_{t+1})    
\end{align} 
\end{subequations}

For the slowed version, we only apply the updates whenever $t\in\mathcal{T}$ and the dynamic regret incurred corresponds to the standard definition used for centralized algorithms; see e.g. \cite{besbes2015non,bedi2018tracking}. We used the ADMM parameters  $\varrho = 1$ and  $\varpi = 0.1$ since they yielded the best performance.

\subsubsection{Tracking via centralized 'slowed' proximal OGD} We consider the proximal OGD algorithm from \cite{dixit2018online}, which proceeds by carrying out the updates in \eqref{ipogd}. As for ADMM, a slowed version is considered and the updates are only carried out at times $t \in \mathcal{T}$. Likewise, the dynamic regret incurred corresponds to the definition used for the centralized case. The step-size parameter was manually set to $\alpha = 0.005$ so as to yield the best performance. 

It is remarked that no comparisons are included with the more general distributed mirror descent algorithm of \cite{shahrampour2018distributed} as it is designed to run on connected graph topologies only. Likewise, performance of SBL-based approaches is not compared since they are generally  not applicable to adversarial settings. 

\subsection{Dynamic regret performance} 
\begin{figure}[t]
	\includegraphics[width=1.09\columnwidth, height=0.25\textheight]{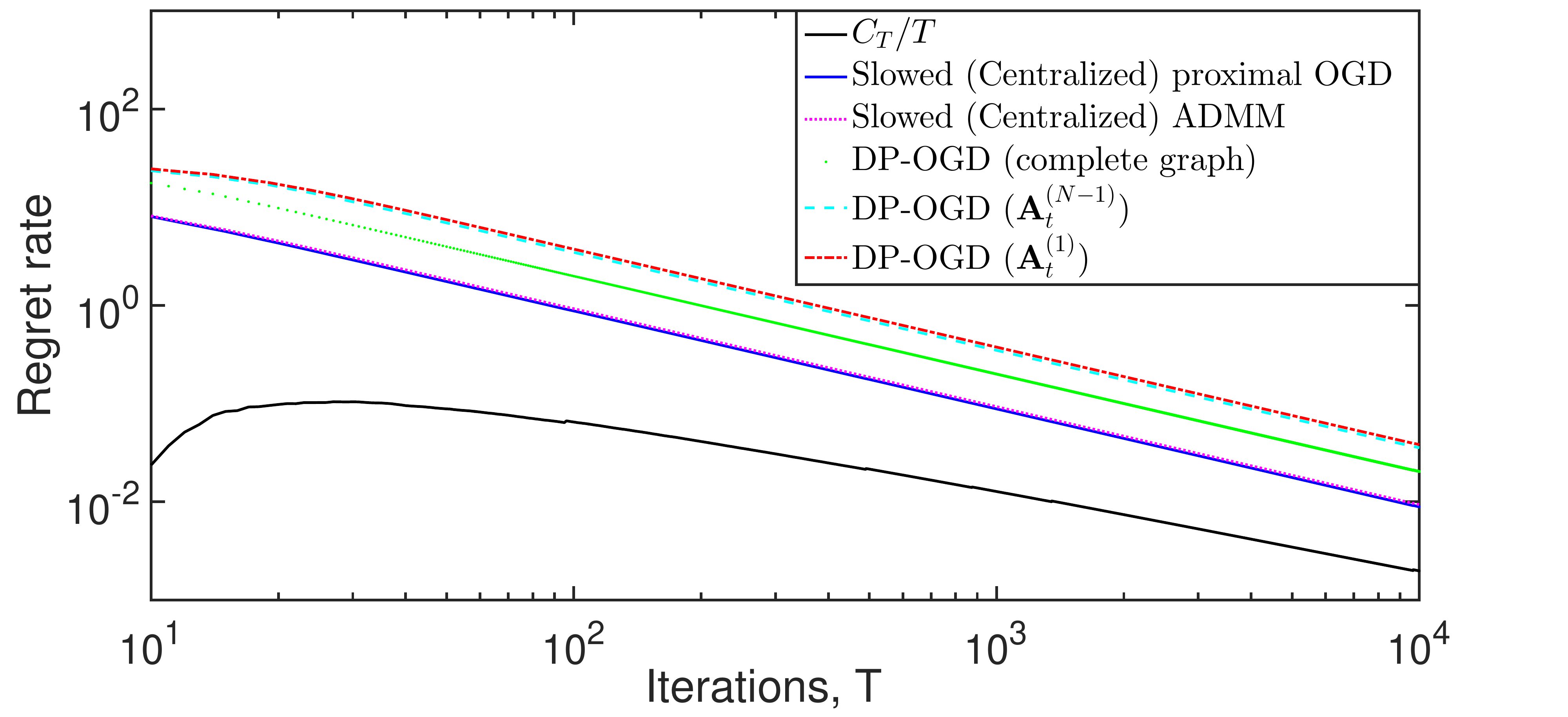}
	\caption{Performance of the distributed and centralized algorithms for $S(k) = 5$.}
	\label{fig1}
\end{figure}
\begin{figure*}[t]
\centering
	\includegraphics[scale=0.3]{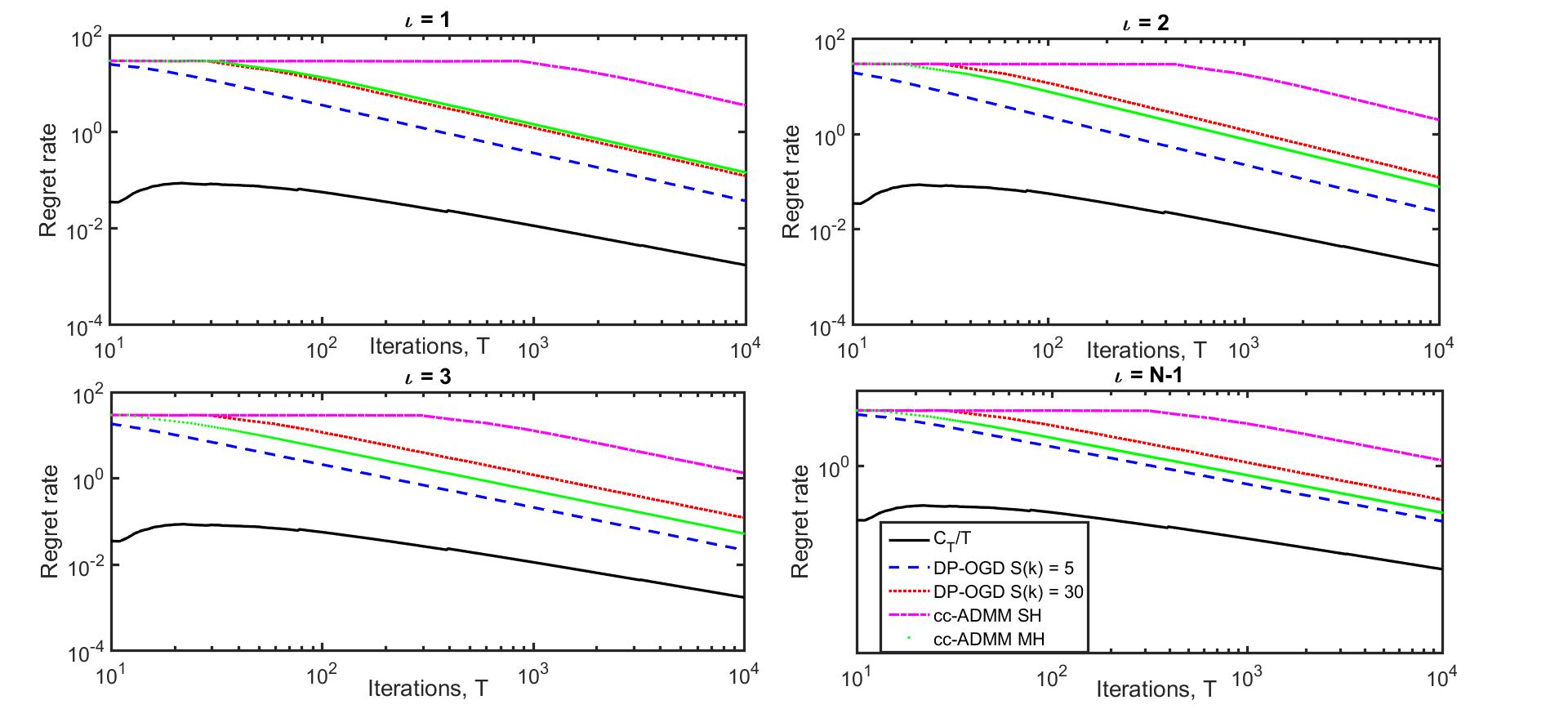}
	\caption{Performance of the distributed and centralized algorithms with matrix $\A_t^{(\iota)}$ for $\iota = 1$, $2$, $3$, and $N-1$. }
	\label{fig2}
\end{figure*}

We begin with considering the case of $S(k) = 5$ for all iterations $k$. \colb{Fig. \ref{fig1}} shows behavior of $\textbf{Reg}_T/T$ against $T$ for the proposed DP-OGD algorithm with complete graph (edge weights $1/N$), single-edge graph matrix $\mathbf{\A}^{(1)}_t$, and the matrix $\mathbf{\A}^{(N-1)}_t$. The regret performance of the centralized 'slowed' schemes is also plotted for comparison. It is remarked that except for the complete graph, it would generally be impossible for all the nodes to exchange updates within $S(k) = 5$ time slots for a network with $N = 100$ nodes. The plot of $C_T/T$ is also shown for comparison with the understanding that the bound in Theorem 1 may include a constant factor. A logarithmic scale is used for plotting so that a regret of $T^\beta$ corresponds to a straight line with slope $\beta-1$. 

It can be seen that the $\textbf{Reg}_T/T$ plots for all the algorithms have the same slope for large $T$. Since the dynamic regret incurred in tracking a time-varying target is always at least $\mathcal{O}(1+C_T)$ even in the centralized case, the plots suggest that the proposed distributed algorithm is close to order optimal. The transient performance of the algorithms is different, and as expected, the 'slowed' centralized ADMM and proximal OGD perform better than the proposed algorithm. Interestingly, the regret incurred when utilizing either $\A_t^{(1)}$ or $\A_t^{(N-1)}$ is almost the same, and not too far from that of the complete graph case. 

As remarked earlier, with $S(k) = 5$, the comparison between distributed and centralized algorithms is not entirely fair, since carrying out the updates for all nodes in just 5 times slots is difficult. Even if the precise implementation details are ignored, since the data $\y_t^i$ is only available at node $i$, the updates require each node to exchange at least some information with every other node, assuming no fusion center is available. Realistically, the time to carry out the updates for the centralized case depends on the network connectivity. For instance, updates can be exchanged within a single time slot when the graph is complete and all nodes can talk to each other for all $t$. More generally, let each pair of nodes be required to exchange a single unit of information at every iteration, and consider the following settings:
\begin{enumerate}
	\item \emph{Single-hop with $\A_{t}^{(\iota)}$}: When only  single hop communication is allowed, each node may transmit the information over edges corresponding to the non-zero entries of $\A_{t}^{(\iota)}$ at time $t$. As a result, the updates must wait till the all information has been exchanged between all pairs of nodes. 
	\item \emph{Multi-hop with $\A_{t}^{(\iota)}$}: When multi-hop communications is allowed, node $i$ sends the update-related information to node $j$ either via a direct path or via a multi-hop path, whichever occurs first. Each edge is allowed to carry a single unit of update-related information each time. Note that such an implementation ignores the fact that distributed implementations of ADMM would generally require at least two exchanges per iteration corresponding to the two primal updates in \eqref{admmup}. 
\end{enumerate}
We refer to these implementations as communication-constrained ADMM (cc-ADMM) for single-hop (SH) and multi-hop (MH) settings. Fig. \ref{fig2} shows the performance of the proposed algorithm for $S(k) = 5$ and $S(k) = 30$, as well as that of the cc-ADMM with SH and MH settings. The four sub-plots correspond to the four different graphs, namely $\A_t^{(1)}$, $\A_t^{(2)}$, $\A_t^{(3)}$, and $\A_t^{(N-1)}$. It can be seen that the proposed algorithm for $S(k) = 5$ is superior in all cases, suggesting that a smaller value of $S(k)$ is preferred in practice. The performance of the centralized algorithm with single-hop communications is poor due to the large delay incurred in the exchange of information among all the nodes. However, the centralized implementation with multi-hop communication comes close to that of the proposed algorithm, though the proposed algorithm uses only single-hop communications.  

Finally, the per-iteration complexity of the ADMM is significantly higher than that of the proposed algorithm. Indeed, since \eqref{admmup} involves minimizing quadratic cost functions, an $n \times n$ system of equations must be solved at every ADMM update, which requires at least $\mathcal{O}(n^3)$ computations per iteration. In contrast, the proposed algorithm requires only $\mathcal{O}(n)$ computations per iteration across all the nodes.

\section{Conclusion}\label{conclusion}
This work considered the problem of tracking the minimum of a time-varying convex constrained optimization problem whose objective function can be expressed \colb{as a} sum of several node-specific costs and a non-differentiable regularizer. The nodes are connected over a time-varying and possibly disconnected graph, thereby complicating the exchange of updates among the nodes. This work proposed a multi-step consensus-based proximal gradient descent algorithm for solving dynamic optimization problems in a distributed manner. Restricted to sampling the loss function sporadically, the number of consensus steps per iteration was carefully selected to yield a sublinear dynamic regret. Dynamic regret of the proposed algorithm is characterized and is shown to be close to that of the centralized tracking algorithm. Numerical tests demonstrate the efficacy of the proposed algorithm on the dynamic sparse recovery problem in wireless sensor networks. 
  		
   \appendices
  \section{Preliminaries}
Before detailing the proofs, we state some of the preliminary relationships and inequalities that will be repeatedly used. 	
	
The proximal operator defined in \eqref{prox} satisfies the following two properties:
\begin{itemize}
	\item Non-expansiveness: for all $\x$, $\y$, 
	\begin{align}
	\norm{\prox_{g_t}^\alpha(\x) - \prox_{g_t}^\alpha(\y)}\leq \norm{\x-\y} \label{ne}
	\end{align}
	\item Relationship with the subgradient: if $\u = \prox_{g_t}^\alpha(\x)$, then it holds that
	\begin{align}
	\x-\u \in \partial g_t(\u)\label{subg}
	\end{align}
	where $\partial g_t(\u)$ is the set of subdifferentials of $g_t$ at $\u$. 
\end{itemize}


  \section{Proof of Lemma \ref{lem1} and Corollary \ref{cor1}}\label{proof_first}
The bound on $\norm{\xb_{k+1}-\xh_k^\star}^2$ is developed by first expanding the squares, using the optimality condition in \eqref{opt}, and the form of the updates in \eqref{izup}-\eqref{ixup}. Subsequently, Cauchy-Schwarz and triangle inequalities are utilized to separate the terms containing $\e_k$ and $\ve_k$. Finally, all the terms containing the gradients $\nabla \fh_k$ are bounded using Assumptions \ref{sm}-\ref{con}. No other assumptions are required at this stage, and therefore Lemma \ref{lem1} is applicable to any inexact proximal OGD setting. 
	
\begin{IEEEproof}[Proof of Lemma \ref{lem1}] The distance between $\bar{\xh}_{k+1}$ and current optimal value $\xh_k^\star$ is given by
\begin{align}
\norm{\xb_{k+1}-\xh_k^\star}^2 =& \norm{\prox_{g_k,\epsilon_k}^\alpha(\zb_k) - \xh_k^\star}^2 \label{first}\\
=& \norm{\prox_{g_k}^\alpha(\zb_k) + \ve_k - \xh_k^\star}^2 \label{useve}\\
=& \norm{\prox_{g_k}^\alpha(\zb_k) - \xh_k^\star}^2 + \norm{\ve_k}^2 \nonumber\\
&+ 2\ip{\ve_k, \prox_{g_k}^\alpha(\zb_k) - \xh_k^\star}  \label{expand1} \\
\leq& \norm{\prox_{g_k}^\alpha(\zb_k) - \xh_k^\star}^2 + \norm{\ve_k}^2 \nonumber\\
&+ 2\norm{\ve_k}\norm{\prox_{g_k}^\alpha(\zb_k) - \xh_k^\star}  \label{second} 
\end{align}
where \eqref{useve} follows from \eqref{ve} and \eqref{second} follows from the applying the Cauchy-Schwartz inequality. Next using the optimality condition in \eqref{opt} and the non-expansiveness property of the $\prox_{g_t}^\alpha$ operator, we obtain
\begin{align}
\norm{\xb_{k+1}-\xh_k^\star}^2 \leq & \norm{\zb_k - (\xh_k^\star - \alpha\nabla \fh_k(\xh_k^\star))}^2 + \norm{\ve_k}^2 \nonumber\\
&+ 2\norm{\ve_k}\norm{\zb_k - (\xh_k^\star - \alpha\nabla \fh_k(\xh_k^\star))}.\label{third}
\end{align}
Next using the update equation in \eqref{izup}, it follows that
\begin{align}  \label{eqst}
 \Vert \bar{\xh}_{k+1}&\!-\!\xh_k^\star \Vert^{2}\\ &\hspace{-8mm}\leq  \Vert (\xb_k-\alpha \nabla \fh_k(\xb_k)) - (\xh_k^\star - \alpha\nabla \fh_k(\xh_k^\star)) - \alpha \e_k \Vert^{2} +\Vert \ve_k \Vert^{2}  \nonumber  \\
&\hspace{-8mm}  +2 \norm{\ve_k}\norm{(\xb_k-\alpha \nabla \fh_k(\xb_k)) - (\xh_k^\star - \alpha \nabla \fh_k(\xh_k^\star)) - \alpha\e_k}  \nonumber \\
&\hspace{-8mm}\leq \Vert (\xb_k-\alpha \nabla \fh_k(\xb_k)) - (\xh_k^\star - \alpha\nabla \fh_k(\xh_k^\star))\Vert^{2} + \Vert \alpha \e_k \Vert^{2}\nonumber\\
&\hspace{-10mm} +2 (\norm{\ve_k} + \norm{\alpha\e_k})\norm{(\xb_k-\alpha \nabla \fh_k(\xb_k)) - (\xh_k^\star - \alpha \nabla \fh_k(\xh_k^\star)) }. \nonumber\\
&\hspace{-5mm} +\Vert \ve_k \Vert^{2} + 2\norm{{\alpha\e_k}}\norm{\ve_k}  \label{tri}   
\end{align}
where \eqref{tri} follows from the Cauchy-Schwarz and triangle inequalities. Next, consider the first summand in \eqref{tri}, which can written as
	\begin{align}\label{fourth}
	&\norm{(\xb_k - \xh_k^\star)  - \alpha(\nabla \fh_k(\xb_k) - \nabla \fh_k(\xh_k^\star))}^2 \nonumber \\
	& =  \norm{\xb_k - \xh_k^\star}^2 + \alpha^2 \Vert \nabla \fh_k(\xb_k) - \nabla \fh_k(\xh_k^\star)\Vert^{2}  \nonumber 
	\\
	& \ \ \ \ -2\alpha \ip{\nabla \fh_k(\xb_k) - \nabla \fh_k(\xh_k^\star) ,\xb_k - \xh_k^\star}
	\end{align}
Since $\{\fh_k^i\}$ are $L$-smooth and $\mu$-convex from Assumptions \ref{sm}-\ref{con}, the same holds for the average function $\fh_k$, allowing us to bound the right-hand side of \eqref{fourth} as
	  \begin{align}	  		 \Vert (&\xb_k - \xh_k^\star)  - \alpha(\nabla \fh_k(\xb_k) - \nabla \fh_k(\xh_k^\star)) \Vert^{2}\nonumber\\	   &\leq  \Vert \xb_k - \xh_k^\star \Vert^{2} + \alpha^{2} L^{2} \Vert \xb_k - \xh_k^\star \Vert^{2} - 2 \alpha \mu \Vert \xb_k - \xh_k^\star \Vert^{2}  \nonumber
 \\
&\leq  \rho^2 \norm{\xb_k - \xh_k^\star}^{2} \label{fifth} 
	  	\end{align}
where, $\rho^2:=(1 + \alpha^2L^2 - 2 \alpha\mu )$. The corresponding bound on the third summand in \eqref{tri} also follows from taking the positive square root in \eqref{fifth}. Applying the bounds in \eqref{fifth} to \eqref{tri} and collecting the terms, we obtain 	  				
\begin{align}
\norm{\bar{\xh}_{k+1}-\xh_k^\star}^{2}  \leq &  \rho^{2} \norm{\xb_k - \xh_k^\star}^{2} + (\norm{\ve_k}+ \norm{\alpha\e_k})^{2} \nonumber\\& + 2 (\norm{\ve_k} + \norm{\alpha\e_k}) \rho \norm{\xb_k-\xh_k^\star}  \\
 = & (\rho \norm{\xb_k-\xh_k^\star} + \norm{\ve_k} + \norm{\alpha\e_k})^{2}   \\
= & (\rho \norm{\xb_k-\xh_k^\star} + \delta_k)^{2}  \label{final}
\end{align}
where $\delta_k:=\norm{\ve_k} + \norm{\alpha\e_k}$. Finally, the required result in Lemma \ref{lem1} follows from taking the positive square root in \eqref{final}. \end{IEEEproof}   

\begin{IEEEproof}[\colb{Proof of Corollary \ref{cor1}}]
The bound on the cumulative sum simply follows from taking sum in \eqref{lemi} over $1\leq k \leq K$, applying triangle inequality to separate the terms containing $\norm{\xh_{k+1}^\star-\xh_k^\star}$ and rearranging. To this end, consider
	 \begin{align}\label{cor_first}
	\!\!\!\!\sum_{k=1}^K\norm{\xb_k- \xh_k^\star } = & \sum_{k=1}^K\norm{\xb_k- \xh_{k-1}^\star + \xh_{k-1}^\star - \xh_k^\star}\nonumber\\
	 \leq & \sum_{k=1}^K \norm{ \bar{\xh}_k-  \xh_{k-1}^\star} + 	\sum_{k=1}^K\norm{\xh_k^\star - \xh_{k-1}^\star}  
	 \end{align}
	 where the second inequality in \eqref{cor_first} uses the triangle inequality. Next utilizing the result of Lemma 1 into \eqref{cor_first}, we obtain
	 \begin{align}\label{cor_second}
	\sum_{k=1}^K\norm{\xb_k- \xh_k^\star } \leq & \sum_{k=1}^K (\rho \norm{\xb_{k-1}- \xh_{k-1}^\star} + \delta_{k-1})  \nonumber \\&+ 	\sum_{k=1}^K\norm{\xh_k^\star - \xh_{k-1}^\star}
\end{align}
Replacing $k-1$ with $k$ in the first summand of \eqref{cor_second} and introducing the $K$-th term on the right-hand side, we obtain
	 \begin{align}\label{cor_third}
	 \sum_{k=1}^K\norm{\xb_k- \xh_k^\star } \leq & \sum_{k=1}^K (\rho \norm{\xb_k- \xh_k^\star} + \delta_k) + 	\sum_{k=1}^K\norm{\xh_k^\star - \xh_{k-1}^\star}\nonumber\\
	&+(\rho \norm{\xb_{0} - \xh_{0}^\star} + \delta_{0}) 
	 \end{align}
If $0 < \rho < 1$ or equivalently if $0< \alpha < 2\mu/L^2$, the first term can be taken to the left to yield 
	\begin{align}
	\sum_{k=1}^K\norm{\xb_k- \xh_k^\star }  \leq & \frac{\rho}{1-\rho}\norm{\xb_{0} - \xh_{0}^\star} + 	\sum_{k=1}^K\frac{\norm{\xh_k^\star - \xh_{k-1}^\star}}{1-\rho}\nonumber\\&  +  \frac{1}{1-\rho}\sum_{k=0}^K \delta_k
	\end{align}
	which is the required result. 
\end{IEEEproof}

\section{Proof of Lemma \ref{cor2}} \label{proofcor2}
The bound in \eqref{cor2a} follows from the update rules in \eqref{zkup}-\eqref{xkup} and the use of the bound in \eqref{wbound}. 

\begin{IEEEproof}[Proof of \eqref{cor2a}] From the definition of $\xb_k$ and the update rule in \eqref{xkup}, we have the following bound
\begin{align}
\sum_{i=1}^{N}\norm{\xb_k- \xh_k^{i} }  &= \sum_{i=1}^{N}\norm{\bigg(\frac{1}{N}\sum_{j=1}^{N}\xh_k^{j}\bigg)- \xh_k^{i} } \label{start1}\\	  		 	
& \hspace{-1.8cm}\leq \frac{1}{N}\sum_{i=1}^{N}\sum_{j=1}^{N}\norm{\xh_k^{i}- \xh_k^{j}}  \label{tri2}\\	
&\hspace{-1.8cm}\leq \frac{1}{N}\sum_{i=1}^{N}\sum_{j=1}^{N}\norm{\prox^{\alpha}_{\gh_{k-1}}(\yh_{k-1}^{i})- \prox^{\alpha}_{\gh_{k-1}}(\yh_{k-1}^{j})}  
\end{align}
where we have used the triangle inequality in \eqref{tri2}. Further use of the triangle inequality and \eqref{ne}, we obtain
 \begin{align}
 \sum_{i=1}^{N}&\norm{\xb_k- \xh_k^{i} }
 \leq  \frac{1}{N}\sum_{i=1}^{N}\sum_{j=1}^{N}\norm{\yh_{k-1}^{i}- \yh_{k-1}^{j}}  \nonumber\\
= & \frac{1}{N}\sum_{i=1}^{N}\sum_{j=1}^{N}\norm{\yh_{k-1}^{i}- \zb_{k-1} + \zb_{k-1} - \yh_{k-1}^{j}} \nonumber \\
 \leq & \frac{1}{N}\sum_{i=1}^{N}\sum_{j=1}^{N}\norm{\yh_{k-1}^{i}- \zb_{k-1}} +  \norm{\zb_{k-1} - \yh_{k-1}^{j}} \nonumber \\
 \leq & 2\sum_{i=1}^{N}\norm{\yh_{k-1}^{i}- \zb_{k-1}}. \label{cor2_last}
\end{align}
Next, consider each summand on the right of \eqref{cor2_last}. Using the update rule in \eqref{zkcon} and the triangle inequality for each $i\in\N$, we have that
\begin{align}
\norm{\yh_{k-1}^{i}- \zb_{k-1}}  = &    \Big\lVert \sum_{j=1}^{N}\left(\colb{Q}^{ij}_{k-1} - \tfrac{1}{N}\right)\zh_{k-1}^{j} \Big\rVert \nonumber\\
 \leq & \sum_{j=1}^{N}\Big\lVert\left(\colb{Q}^{ij}_{k-1} - \tfrac{1}{N}\right)\zh_{k-1}^{j} \Big\rVert \nonumber\\
\leq & \sum_{j=1}^{N}\abs{\colb{Q}^{ij}_{k-1} - \frac{1}{N}}\norm{\zh_{k-1}^{j} }
 \end{align}
It is now possible to utilize the bound in \eqref{wbound} (see Assumptions \ref{bcon}-\ref{wt}) as follows. 
  \begin{align}\label{last}
  \norm{\yh_{k-1}^{i}- \zb_{k-1}} \leq &  \sum_{j=1}^{N} \Gamma \gamma^{S(k-1)-1} \norm{\zh_k^{j} }\nonumber \\
 \leq & \Gamma \gamma^{S(k-1)-1}  \sum_{j=1}^{N}\norm{\zh_{k-1}^{j} } 
 \end{align}
Since the right-hand side is the same for all $i\in\N$, substituting \eqref{last} into \eqref{cor2_last}, we obtain the desired result where the bound includes an additional factor of $2N$.  
\end{IEEEproof}
 
\begin{IEEEproof}[Proof of \eqref{cor2b}]
Next, we establish that the $\sum_{i=1}^{N}\norm{ \zh_{k-1}^{i}}$ is upper bounded by an affine function of $k$. The result follows in a straightforward manner from the fact that every update entails adding bounded terms. Specifically, taking norm and applying the triangle inequality to \eqref{zkup}, we obtain
\begin{align}\label{gradient_bound}
\norm{\zh_{k-1}^i} \leq & \norm{\xh_{k-1}^i} + \alpha \norm{\nabla \fh_{k-1}^i(\xh_{k-1}^i)} \\
\leq &  \norm{\xh_{k-1}^i} + \alpha M \label{gbound1}
\end{align}
where \eqref{gbound1} follows from Assumption \ref{lip}. Taking summation over $i\in\N$, we obtain
\begin{align}
\sum_{i=1}^{N}\norm {\zh_{k-1}^{i}} \leq & \sum_{i=1}^{N}\norm{\xh_{k-1}^{i}} + \alpha NM.    \label{telescopic_sum}
\end{align}
In order to develop an upper bound on the first term in \eqref{telescopic_sum}, we use the property \eqref{subg} to obtain
\begin{align}
\xh_{k-1}^i = \yh_{k-2}^i - \alpha \partial \gh_{k-2}(\xh_{k-1}^{i})
\end{align}
Again taking norm on both sides and using triangle inequality, we get
\begin{align}
\norm{ \xh_{k-1}^{i}} \leq & \norm{\yh_{k-2}^{i}} + \alpha \norm{\partial \gh_{k-2}(\xh_{k-1}^{i})} \\
&\leq  \norm{\yh_{k-2}^{i}} + \alpha M
\end{align}
where the bound on the subgradient follows from \eqref{lipg}. Observe further that $\yh_{k-2}^i$ is a convex combination of $\{\zh_{k-2}^i\}_{i=1}^N$ since $\yh_{k-2}^i = \sum_{j=1}^{N}\colb{Q}^{ij}_{k-2}\zh_{k-2}^j $ and $\colb{\Q}_k$  is doubly stochastic as it is a product of doubly stochastic matrices. Therefore from the convexity of the norm, it follows that $\norm{\yh_{k-2}^i} \leq \sum_{j=1}^N \colb{Q}_k^{ij}\norm{\zh_{k-2}^i}$, implying that
\begin{align}
\sum_{i=1}^N\norm{ \xh_{k-1}^i} \leq & \sum_{i=1}^{N}\sum_{j=1}^N\colb{Q}_k^{ij}\norm{\zh_{k-2}^j} + \alpha N M  \nonumber \\
\leq & \sum_{i=1}^{N}\norm{\zh_{k-2}^{i}} + \alpha N  M  \label{convex_combination}
\end{align} 
Substituting \eqref{telescopic_sum} into \eqref{convex_combination}, we obtain
\begin{align}
\sum_{i=1}^{N}\norm{ \zh_{k-1}^{i}} \leq & \sum_{i=1}^{N}\norm{\zh_{k-2}^{i}} +  2\alpha N  M  \\
\leq & \sum_{i=1}^{N}\norm{\zh_{0}^{i}} +  2\alpha N  M  (k-1) 
\end{align}
which is the required result.
\end{IEEEproof}
  		   		 

\section{\colb{Proof of Lemma \ref{lemdel}}} \label{lemma_2_proof}
The proof is split into two parts, corresponding to obtaining bounds on $\norm{\e_k}$ and $\norm{\ve_k}$. 

\begin{IEEEproof}[Bound on $\norm{\e_k}$] Recall that the gradient error $\e_k$ is defined as
\begin{align}
\e_k &= \frac{1}{N}\sum_{i=1}^{N}(\nabla \fh_k^{i}(\xh_k^{i}) - \nabla \fh_k^{i}(\xb_k) )  .
\end{align}
Taking norm on both sides, using triangle inequality, and Assumption \ref{sm}, we obtain 
\begin{align}
\norm{ \e_k}\leq &    \frac{1}{N}\sum_{i=1}^{N}\norm{(\nabla \fh_k^{i}(\xh_k^{i}) - \nabla \fh_k^{i}(\xb_k) )}  \nonumber\\
\leq &    \frac{1}{N}\sum_{i=1}^{N}L\norm{\xh_k^{i} - \xb_k} \end{align}
Next, using the upper bound established in Lemma \ref{cor2}, the required bound becomes
\begin{align}\label{delta_first}
\norm{ \e_k}&\leq  2 L  \Gamma \gamma^{S(k-1)-1}  \sum_{j=1}^{N}\norm{\zh_{k-1}^{j} }\nonumber\\
&\leq2 L \Gamma \gamma^{S(k-1)-1}  \left[\sum_{i=1}^N\norm{\zh_{0}^{i}} +  2\alpha N M (k-1)\right].
\end{align} 
\end{IEEEproof}

\begin{IEEEproof}[Bound for $\norm{\ve_k}$] Using the definition of $\ve_k$ in \eqref{ve} and the triangle inequality, we obtain
\begin{align}
\norm{\ve_k} &= \norm{\frac{1}{N}\sum_{i=1}^N\prox_{\gh_k}^\alpha(\y_k^i)-\prox_{\gh_k}^{\alpha}(\zb_k)} \\
&\leq \frac{1}{N}\sum_{i=1}^N \norm{\prox_{\gh_k}^\alpha(\y_k^i)-\prox_{\gh_k}^{\alpha}(\zb_k)} \\
&\leq \frac{1}{N}\sum_{i=1}^N \norm{\y_k^i-\zb_k} \label{neprox2}
\end{align}
where \eqref{neprox2} follows from \eqref{ne}. Therefore, using the bound in \eqref{last} we have that
\begin{align}
\norm{\ve_k} &\leq \Gamma \gamma^{S(k)-1}  \sum_{j=1}^{N}\norm{\zh_k^{j} } \\
&\leq \Gamma \gamma^{S(k)-1}\left[\sum_{i=1}^N \norm{\zh_0^i} + 2\alpha NMk\right]
\end{align}
The required bound on $\delta_k$ follows from combining the two parts of the proof. 
\end{IEEEproof}
  

\section{Proof of Theorem 1} \label{Proof_thm_1}For the sake of compactness, let $\ell_t^i:= f_t^i + g_t$. Recall that $\x_t^i = \x_{\t}^{\colb{i}} = \x_{t_k}^i$ for $k$ such that $t_k \leq t < t_{k+1}$. Using the first order convexity condition, we have that
\begin{align}
 \Reg=& \frac{1}{N^2} \sum_{t=1}^T\sum_{i=1}^{N} \sum_{j=1}^{N}\left(\ell_t^j(\x_t^i) - \ell_t^j(\x_t^\star)\right) \nonumber \\
\leq & \frac{1}{N^2}\sum_{t=1}^{T}\sum_{i=1}^{N}\sum_{j=1}^{N}\ip{{\nabla} \ell_t^{j}(\x^{i}_{\t}), \x_{\t}^{i}-\x_t^\star}\label{first_thm}
\end{align}
Using Cauchy Schwartz inequality and gradient boundedness on the right hand side of \eqref{first_thm}, we obtain 
 \begin{align} 
\Reg\leq & \frac{1}{N^2}\sum_{t=1}^{T}\sum_{i=1}^{N}\sum_{j=1}^{N} \norm{{\nabla} \ell_t^{j}(\x^{i}_{\t})}  \norm{\x^{i}_{\t}- \x_t^\star}\nonumber\\
& \leq \frac{2M}{N}\sum_{t=1}^{T}\sum_{i=1}^{N}\norm{\x^{i}_{\t}- \x_{t}^\star}\\
&\leq \frac{2M}{N}\sum_{t=1}^{T} \sum_{i=1}^{N}(\norm{\x^i_{\t}- \x_{\t}^\star} + \norm{\x_{\t}^\star - \x_t^\star}) \nonumber\\
&\hspace{-1.2cm} =\frac{2M}{N}\sum_{t=1}^{T} \sum_{i=1}^{N} \norm{\x^i_{\t}- \x_{\t}^\star} + 2M\sum_{t=1}^T \norm{\x_{\t}^\star - \x_t^\star}.\label{regret_bound}
\end{align}
It can be seen that the dynamic regret is a sum of two components, namely terms depending on the optimality gap and on the path length. Of these, the first component only contains the optimality gap evaluated at times $\t \in \T$ since the agents takes actions only at those times. The second term in the \eqref{regret_bound} also contains the error incurred due to sampling the objective function intermittently and depends on both, the path length and the sequence $\{S(k)\}$ as will subsequently be shown. 
  		  	 
First consider the optimality gap in \eqref{regret_bound}, which can be bounded by first converting the indexing to $k$ and then adding and subtracting the network averaged iterates $\xb_k$. For any $t\geq 1$, it holds that 
\begin{align}
\sum_{t=1}^T&\sum_{i=1}^{N}\norm{\x^i_{\t}- \x_{\t}^\star} = \sum_{k=1}^K\sum_{i=1}^N (S(k)+2)\norm{\xh_k^i-\xh_k^\star} \\
&\leq \sum_{k=1}^K\sum_{i=1}^N (S(k)+2)(\norm{\xh_k^i-\xb_k} + \norm{\xb_k - \xh_k^\star})\\
&\leq (S(K)+2) \sum_{k=1}^K\sum_{i=1}^N (\norm{\xh_k^i-\xb_k} + \norm{\xb_k - \xh_k^\star})\label{tri3}
\end{align}
where \eqref{tri3} follows from the triangle inequality and the fact that $S(k)$ is non-decreasing in $k$ so that $S(k)\leq S(K)$ for all $k\leq K$. The first term in \eqref{tri3} can be bounded in from Lemma \ref{cor2} while the second term can be bounded from Corollary \ref{cor1} to yield
\begin{align}
\sum_{t=1}^T&\sum_{i=1}^{N}\norm{\x^i_{\t}- \x_{\t}^\star} \nonumber\\
&\leq \sum_{t=1}^{T} 2 \Gamma \gamma^{S(k-1)-1} N \bigg[\sum_{i=1}^{N}\norm{\zh_{0}^{i}} + 2\alpha N M (k-1) \bigg ] \nonumber\\
&+ \frac{N(S(K)+2)}{1-\rho}\left( \rho\norm{\xb_0 - \xh_0^\star} + C_T+E_T\right)\label{gap1}\\
&=\mathcal{O}(R_T(1+E_T+C_T))
\end{align}
where \eqref{gap1} also uses the bounds in Lemma \ref{lemdel}. The second term in \eqref{regret_bound} can be bounded by observing that
\begin{align}
\sum_{t=1}^T\norm{\x_{\t}^\star - \x_t^\star} &\leq \sum_{t=1}^T\sum_{\tau = \t}^{t-1} \norm{\x_{\tau+1}^\star - \x_\tau^\star} \\
&\leq (S(K)+2)C_T = \mathcal{O}(R_TC_T)\label{err2bnd}
\end{align}
Combining the bounds in \eqref{gap1} and \eqref{err2bnd}, and substituting $S(K) + 2 = R_T$, we obtain the desired regret bound as
\begin{align}
\Reg \leq \mathcal{O}(R_T(1+E_T + C_T)).
\end{align}

\footnotesize
\bibliographystyle{IEEEtran} 
\bibliography{IEEEabrv,ref}
	\begin{IEEEbiography}
								{Rishabh Dixit} 
								received  his B.Tech degree in Electrical Engineering from Indian Institute of Technology, Kanpur in 2015. Currently, he is pursuing the Ph.D. degree with the Department of Electrical \& Computer Engineering at Rutgers University, NJ, USA. 
								
							 His research interests include algorithms development using convex optimization, dynamic optimization, and and machine learning.  
							\end{IEEEbiography}
							
	\begin{IEEEbiography}
								{Amrit Singh Bedi} (S'16)
								received  the Diploma degree in electronics and communication engineering (ECE) from Guru Nanak Dev Polytechnic, Ludhiana, India, in 2009, the B.Tech. degree in ECE from the Rayat and Bahra Institute of Engineering and Bio-Technology, Kharar, India, in 2012, and the M.Tech. degree in electrical engineering (EE) from IIT Kanpur, in 2016, where he is currently pursuing the Ph.D. degree with the EE Department. 
								
								He is currently involved in developing stochastic optimization algorithms for supervised learning. His research interests include distributed stochastic optimization for networks, time-varying optimization, and machine learning.  
							\end{IEEEbiography}

								\begin{IEEEbiography}
								{Ketan Rajawat}
								(S'06-M'12) received his B.Tech and M.Tech degrees in Electrical Engineering from the Indian Institute of Technology (IIT) Kanpur, India, in 2007, and his Ph.D. degree in Electrical and Computer Engineering from the University of Minnesota, USA, in 2012. Since 2012, he has been an Assistant Professor at the Department of Electrical Engineering, IIT Kanpur. His current research focuses on distributed algorithms, online learning, and control of networked systems. He is currently an Associate Editor of the IEEE Communications Letters. He is a recipient of the Indian National Science Academy (INSA) Medal for Young Scientists in 2018 and the P. K. Kelkar Young Faculty Research Fellowship in 2018. 
								
							\end{IEEEbiography}
\end{document}